\documentclass[oneside,english]{amsart}
\usepackage[T1]{fontenc}
\usepackage[latin9]{inputenc}
\usepackage{float}
\usepackage{amstext}
\usepackage{amsthm}
\usepackage{amssymb}
\usepackage{graphicx}

\makeatletter

\providecommand{\tabularnewline}{\\}

\numberwithin{equation}{section}
\numberwithin{figure}{section}
\theoremstyle{plain}
\newtheorem{thm}{\protect\theoremname}[section]
\theoremstyle{definition}
\newtheorem{defn}[thm]{\protect\definitionname}
\theoremstyle{plain}
\newtheorem{question}[thm]{\protect\questionname}
\theoremstyle{remark}
\newtheorem{rem}[thm]{\protect\remarkname}

\usepackage{yhmath}
\usepackage{hyperref}
\usepackage{breqn}

\makeatother

\usepackage{babel}
\providecommand{\definitionname}{Definition}
\providecommand{\questionname}{Question}
\providecommand{\remarkname}{Remark}
\providecommand{\theoremname}{Theorem}

\begin{document}
\global\long\def\Rd{\mathbb{R}^{d}}

\global\long\def\R{\mathbb{R}}

\global\long\def\Z{\mathbb{Z}}

\global\long\def\N{\mathbb{N}}

\global\long\def\parenth#1{\left(#1\right)}

\global\long\def\radii{\textup{radii}}

\global\long\def\centers{\textup{centers}}

\global\long\def\codecenterlabel{\textup{center}\lambda}

\global\long\def\codeneighborlabels{\lambda}

\global\long\def\neighbors{\textup{neighbors}}

\global\long\def\cal#1{\mathcal{#1}}

\global\long\def\set#1#2{\left\{  \vphantom{{#1#2}}#1\right|\left.\vphantom{{#1#2}}#2\right\}  }

\global\long\def\mickeymouse#1#2#3{\widehat{#2#1#3}}

\global\long\def\teddybear#1#2#3{\widehat{#1#2/#3}}

\global\long\def\tripod#1#2{\widehat{#1/#2}}

\global\long\def\necklace#1#2#3{#1#2{:}#3}

\global\long\def\pooh#1#2#3{P{#2\atop #1}\mathrel{:}#3}

\global\long\def\tripods{\textup{tripods}}

\global\long\def\necklaces{\textup{necklaces}}

\global\long\def\teddybears{\textup{teddybears}}

\global\long\def\angle{\textup{ang}}

\global\long\def\dihed{\textup{dihed}}

\global\long\def\solid{\textup{solid}}

\newcommand{\miek}{Miek Messerschmidt}

\newcommand{\miekemail}{mmesserschmidt@gmail.com}

\newcommand{\witsaddress}{
School of Mathematics, 
University of the Witwatersrand, 
Private Bag~3, 
Wits~2050, 
South Africa}

\newcommand{\subjectcodesprimary}{05B40, }

\newcommand{\subjectcodessecondary}{52C17, }

\newcommand{\keywordstext}{raspberries, sphere packing, compact sphere packings}

\title[Raspberries with two sizes of berry]{The raspberries in three dimensions\\ with at most two sizes of berry}

\author{\miek}
\address{\miek,\ \witsaddress}
\email{\miekemail}

\begin{abstract}
In three dimensional Euclidean space, a raspberry is defined to be
an arrangement of spheres with pairwise disjoint interiors, so that
all spheres are tangent to a central unit sphere and is such that
the contact graph of the non-central spheres triangulate the central
sphere.

We discuss the relevance of these structures in other work. We present
a catalog of all configurations of radii that permit the formation
of raspberries that have at most two sizes of non-central spheres.
Throughout we discuss the construction of this catalog. \end{abstract}

\subjclass[2020]{Primary: \subjectcodesprimary. Secondary: \subjectcodessecondary}

\maketitle

\section{Introduction}

Sphere packing has a long history with contributions from many prominent
mathematicians. In this subfield of mathematics, the so-called compact
sphere packings\footnote{A compact sphere packing of $d$-dimensional Euclidean space is such
that its contact hypergraph triangulates the entire ambient space.} (cf. \cite{Fernique2021,FerniqueTwoSpheres2021,FerniqueThreeSizes,Kennedy2006,Messerschmidt2020,Messerschmidt2d,MesserschmidtKikianty2024})
are of interest and have relevance, for example, in materials science
where such structures have been observed `in the wild' (cf. \cite{ChemistryPaik}
and \cite{ChemistryFernique}). 

\begin{figure}
\includegraphics[totalheight=4cm]{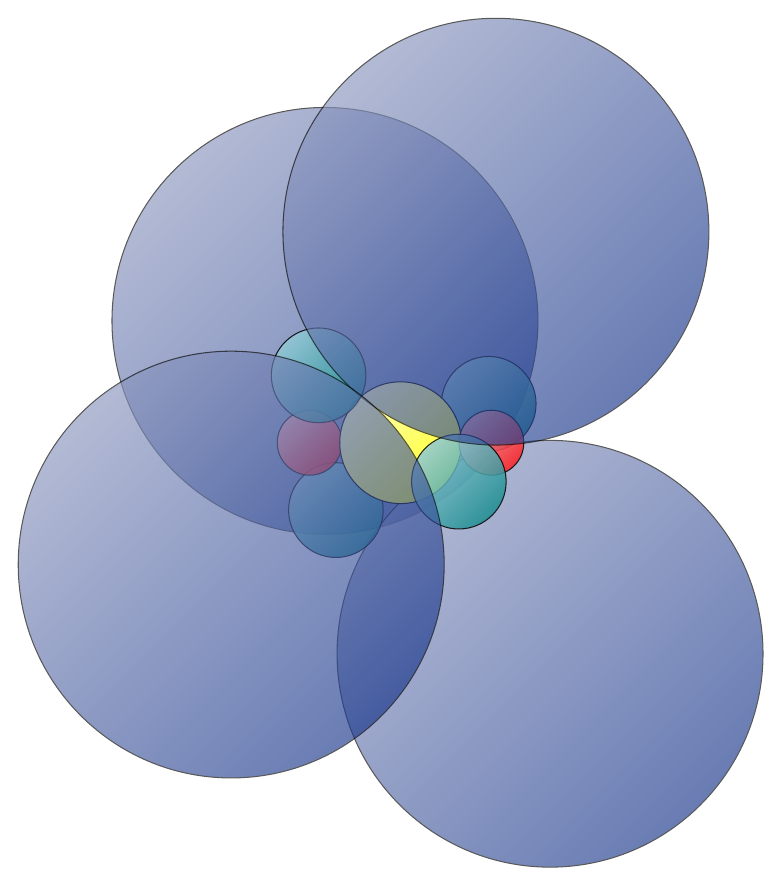}\includegraphics[totalheight=4cm]{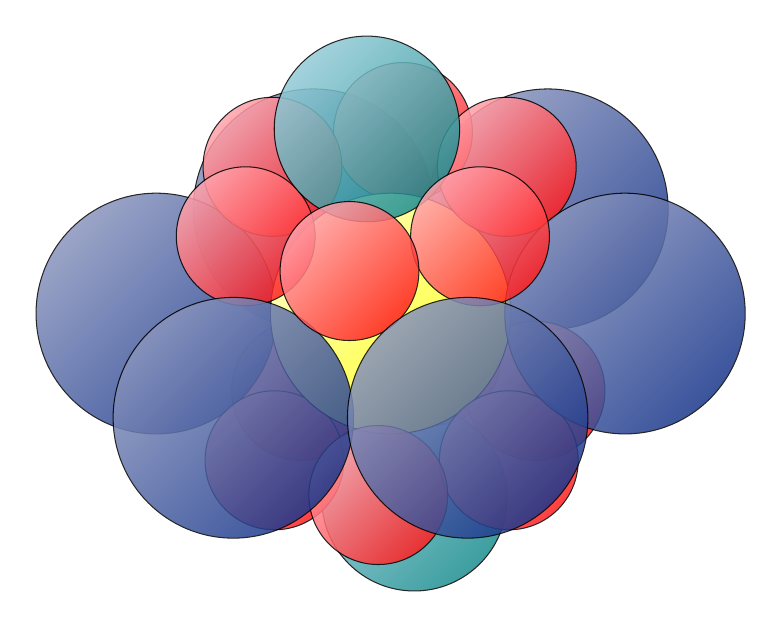}

\caption{Examples of raspberries with three sizes of berry. \label{fig:example-raspberry}}
\end{figure}

In this paper we will focus on what we call \emph{raspberries}\footnote{We are not aware of any previously established terminology. The name
`raspberry' is chosen for these structures' vague resemblance to the
fruit. } which are important structures that occur in compact sphere packings
of three-dimensional Euclidean space, cf. \cite{FerniqueTwoSpheres2021,FerniqueThreeSizes,MesserschmidtKikianty2024},
but turn out to be interesting in their own right. Raspberries are
the three-dimensional analogue of `flowers' or `coronas' in two dimensions
which occur in compact disc packings of the Euclidean plane, (cf.
\cite{Kennedy2006,Fernique2021,MathewsZymaris}). Figure~\ref{fig:example-raspberry}
displays two examples of raspberries. 

We provide a formal definition:
\begin{defn}
\label{def:raspberry}In three dimensional Euclidean space we define
a \emph{raspberry }to be an arrangement of finitely many spheres that
satisfy the following three conditions:\vspace{2mm}
\begin{enumerate}
\item All spheres in the arrangement have pairwise disjoint interiors. \vspace{2mm}
\item A unit sphere, centered at the origin, occurs in the arrangement,
and all other spheres are tangent to this central unit sphere. The
central sphere is called the \emph{pit }of the raspberry.\emph{ }The
non-central spheres are called the \emph{berries} of the raspberry.\vspace{2mm}
\item \label{enu:no-hemispheres-or-lunes}The contact graph of the berries
is formed by taking the centers of the berries as vertices, and if
two berries are mutually tangent, connecting their centers by an edge
(represented by a line segment). Through central projection onto the
pit, the contact graph of the berries is required to triangulate the
pit with spherical triangles. We call this triangulation of the pit
the \emph{berry triangulation}. We note that hemispheres and lunes
are not considered spherical triangles. Therefore there exists no
closed half space that contains all the vertices of the contact graph
of the berries.
\end{enumerate}
\end{defn}

It is possible to construct raspberries by guessing an appropriate
combinatorial structure and solving for the sizes of the berries (see
Figure \ref{fig:example-raspberry}). However, to our knowledge, little
is known of raspberries in general. In particular, a currently existing
blind spot for raspberries is the lack of an analogue to Mathews and
Zymaris' recent generalization of Descartes' Theorem for flowers \cite{MathewsZymaris}.

\hspace{5mm}

In this paper we construct \emph{all} configurations of radii that
permit the formation of raspberries that have at most two sizes of
berry. The main purpose of this work, apart from providing a catalog
akin to those in \cite{Kennedy2006,Fernique2021,FerniqueTwoSpheres2021,FerniqueThreeSizes},
is to additionally investigate, in a concrete manner, the geometric
properties of  raspberries in the simplest nontrivial case. To elaborate,
our particular interest in the geometric properties of raspberries
stems from \cite[Corollary~6.9]{MesserschmidtKikianty2024} where
it is proved that in compact sphere packings of $d$-dimensional Euclidean
space with at most $n$ sizes of sphere, there are at most finitely
many configurations of radii of spheres allowed, granted that a so-called
\emph{heteroperturbative hypothesis }is additionally assumed. A family
of triangulations of the unit $(d-1)$-sphere is said to be \emph{heteroperturbative}
if, upon perturbation of any triangulation in the family that preserves
the combinatorics, but changes the lengths of some of the edges of
the triangulation while remaining inside the family, the perturbed
triangulation must have at least one edge longer and at least one
edge shorter than the corresponding edges in the original triangulation
(cf. \cite[Section~2.4]{MesserschmidtKikianty2024}). This raises
the following still-open question:
\begin{question}
\label{que:berry-triangulations-heteroperturbative}Is the family
of all berry triangulations heteroperturbative?
\end{question}

Given the definition of a raspberry and its berry triangulation, this
question has a clear generalization to arbitrary dimensions. We remark
that in the Euclidean plane, it is easily seen that the family of
all triangulations of the unit $1$-sphere (by circular arcs) is heteroperturbative
\cite[Example~2.3]{MesserschmidtKikianty2024}. In higher dimensions,
the answer is not clear and thus provides a motivation for the current
paper. In contrast to triangulations of the $1$-sphere, it is not
true of the family of \emph{all} triangulations of the unit $2$-sphere
in three-dimensional Euclidean space is heteroperturbative \cite[Example~2.4]{MesserschmidtKikianty2024}.
This is not to say that this question is moot in higher dimensions,
as a result of Winter shows that the family of triangulations of the
unit $(d-1)$-sphere that arise through central projection of inscribed
\emph{convex} (simplicial) polytopes with the center of the inscribing
sphere in the polytope's interior, is heteroperturbative \cite[Corollary~4.13]{WinterPolytopes}.
Unfortunately not all berry triangulations arise in this way from
a convex polytope \cite[Figure~6]{MesserschmidtKikianty2024}. Still,
the results in Sections~\ref{sec:Computing-all-raspberries} and~\ref{sec:thirty-nonflexible}
of the current paper partially answers Question~\ref{que:berry-triangulations-heteroperturbative}
affirmatively in the case of berry triangulations arising from raspberries
with at most two sizes. 

\hspace{5mm}

We briefly describe the contents of this paper. Section~\ref{sec:Preliminaries}
provides some preliminary definitions and notation regarding standard
angles that can occur in a raspberry. We mostly follow the notation
and terminology from \cite{FerniqueTwoSpheres2021}. Importantly,
we define what we mean by a canonical labeling of spheres in a raspberry,
by a (01-- and 02--) necklace, and by a (01-- and 02--) necklace
code. In Definition~\ref{def:constraints}, we explicitly state constraints
that raspberries must satisfy. 

Section~\ref{sec:necklace-polynomials} describes the procedure followed
to compute multivariate necklace polynomials (in indeterminates representing
berry radii) from necklace codes. If a particular necklace code is
determined by a necklace in a raspberry, the configuration of radii
of the berries in the raspberry must, by construction, necessarily
lie on the algebraic variety determined by the corresponding necklace
polynomial. Practically, necklace polynomials are computed from necklace
codes through performing some precomputation and two subsequent Gr\"obner
basis computations so as to obtain a polynomial in indeterminates
representing only the berry radii. 

In Section~\ref{sec:Computing-all-raspberries} we describe the procedure
by which all configurations of radii that permit the formation of
raspberries with at most two sizes of berries are computed. As a first
step, in Section~\ref{subsec:The-01-necklace-codes-and-polynomials},
through placing bounds on dihedral angles that can possibly occur
in a 01-necklace, we compute all 01-necklace codes and corresponding
necklace polynomials that can occur in a raspberry with at most two
sizes of berry. These data are presented in Table~\ref{tab:01-necklaces-and-polynomials}.
Before moving on, we present all raspberries with only one size of
berry in Section~\ref{subsec:raspberries-one-size} and all flexible
raspberries (for which berry radii are not uniquely determined) with
at most two sizes of berry in Section~\ref{subsec:raspberries-flex}.
In Section~\ref{subsec:Bounds-on-the-02-dihed-angles}, we determine
bounds on the dihedral angles in a 02-necklace in a hypothetical raspberry,
under the assumption that a 01-necklace with necklace code from Table~\ref{tab:01-necklaces-and-polynomials}
is present in the hypothetical raspberry. These bounds are presented
in Table~\ref{tab:02-teddy-bounds} and by exploiting them in Section~\ref{subsec:computing-pairs}
along with the necklace dihedral angle constraint from Definition~\ref{def:constraints},
using interval arithmetic, we determine 596 pairs of 01-- and 02-necklace
codes that can possibly simultaneously solve the necklace dihedral
angle constraints from Definition~\ref{def:constraints} in the possible
berry radii $r_{1}$ and $r_{2}$. In Sections~\ref{subsec:excluding-necklace-code-pairs}
and Sections~\ref{subsec:analyzing -numerical-computations} we carefully
winnow this number down to 30 pairs of radii $r_{1}$ and $r_{2}$
that can appear as radii in a raspberry. Finally we present a catalog
of 30 raspberries with these sizes of berry in Section~\ref{sec:thirty-nonflexible}.

\hspace{5mm}

Before proceeding, we make a note on the computational tools employed.
We employ Python and the Python libraries \texttt{sympy} \cite{sympy}
and \texttt{mpmath} \cite{mpmath} to perform symbolic computations,
numerical solving, arbitrary precision (interval) arithmetic, and
real root finding of univariate polynomials. To simplify systems of
multivariate polynomial equations, we employ elimination theory (see
for example \cite[Section~2.3]{AdamsLoustaunau}) through computation
of Gr\"obner bases. For this purpose, we employ both \texttt{Singular}
\cite{SINGULAR} and \texttt{sympy}. 

We provide two structured datasets in \texttt{JSON} format that contain
the results of computations that we will reference throughout the
paper and that will permit the reader to independently verify the
claims made \cite{raspberrydata}. 

\section{Preliminaries\label{sec:Preliminaries}}

Let $R$ be a raspberry. We define the canonical labeling of the raspberry
$R$ as follows. The pit $P$ of $R$ is always labeled with the label
`$0$' and we set $r_{0}:=1$. Let $n\in\N$, with $0<r_{1}<r_{2}<\ldots<r_{n}$
denoting the all the values of the radii of all the berries that occur
in the raspberry $R$. Each berry with radius $r_{j}$ then is labeled
by the corresponding index $j$. 

Let $A,B,C$ be mutually tangent spheres in $R$ with respective labels
$a,b$ and $c$. Using the planar law of cosines, the angle, denoted
$\mickeymouse cab$, formed at the center of $C$ by the rays respectively
through the centers of $A$ and $B$ is satisfies 
\[
\cos(\mickeymouse cab)=\frac{(r_{c}+r_{a})^{2}+(r_{c}+r_{b})^{2}-(r_{a}+r_{b})^{2}}{2(r_{c}+r_{a})(r_{c}+r_{b})}.
\]

Let $H,A,B$ be mutually tangent berries in $R$ with respective labels
$h,a,b$. Using the spherical law of cosines, the dihedral angle,
denoted $\teddybear 0h{ab}$, that is formed by the planes respectively
through the centers of the spheres $P,H,A$ and $P,H,B$ satisfies
\[
\cos(\teddybear 0h{ab)}=\frac{\cos(\mickeymouse 0ab)-\cos(\mickeymouse 0ha)\cos(\mickeymouse 0hb)}{\sin(\mickeymouse 0ha)\sin(\mickeymouse 0hb)}.
\]
By Girard's Theorem, the area, denoted $\tripod 0{hab}$, formed by
the spherical triangle on the pit $P$ by the points of contact on
$P$ by $H,A,B$ satisfies
\[
\tripod 0{hab}=\teddybear 0h{ab}+\teddybear 0a{bh}+\teddybear 0b{ah}-\pi.
\]

\begin{figure}
\includegraphics[totalheight=4cm]{example_raspberry02}\includegraphics[totalheight=4cm]{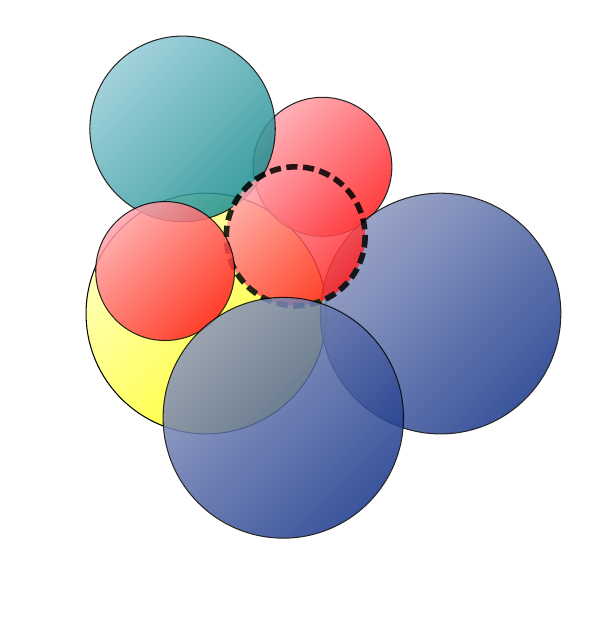}

\caption{A necklace with necklace code $\protect\necklace 01{33121}$ from
a raspberry with three sizes of berry. \label{fig:example-necklace}}
\end{figure}
Each berry $H$ in a raspberry determines a \emph{necklace,} $N$,
meaning a cycle of spheres of some length $m\in\N$, $B_{0},B_{1},\ldots,B_{m-1}$
so that, for all $j\in\Z/m\Z$, the sphere $B_{j}$ is mutually tangent
to the spheres $B_{j+1}$, $H$, and the raspberry's pit. The sphere
$H$ is called the \emph{head }of the necklace $N$, and we call the
spheres $B_{0},B_{1},\ldots,B_{n-1}$ the \emph{beads }of the necklace
$N$ (see Figure~\ref{fig:example-necklace}). With $h$ the label
of the head $H$, and $b_{0},b_{1},\ldots,b_{n-1}$ the labels of
the beads of $N$, we define the \emph{necklace code }of the necklace
$N$ as the symbol $\necklace 0h{b_{0}b_{1}\ldots b_{n-1}}$. We refer
to such a necklace (code) as a $0h$-necklace (code). We regard necklace
codes that differ only in a reflection or rotation of the cycle of
bead labels as identical. 

It is clear that every raspberry necessarily satisfies the following
constraints which we name explicitly in the next definition for later
referral.
\begin{defn}
\label{def:constraints} Let $R$ be a canonically labeled raspberry.\vspace{2mm}
\begin{enumerate}
\item \textbf{The necklace dihedral angle constraint.} For each necklace
$N$ in the raspberry $R$, the sum of the dihedral angles in the
necklace equals $2\pi$. Explicitly, if $\necklace 0h{b_{0}b_{1}\ldots b_{n-1}}$
is the necklace code for the necklace $N$ in the raspberry $R$,
the following equation is satisfied 
\[
\sum_{j\in\Z/m\Z}\teddybear 0h{b_{j}b_{j+1}}=2\pi.
\]
 \vspace{2mm}
\item \textbf{The pit area constraint.} The total area of the triangles
in the berry triangulation equals the area of the pit of the raspberry
$R$, i.e., $4\pi$. For a necklace $N$ in the raspberry $R$, let
$\necklace 0h{b_{0}b_{1}\ldots b_{n-1}}$ denote the necklace code
for $N$ and define 
\[
\text{Area}(N):=\sum_{j\in\Z/m\Z}\tripod 0{hb_{j}b_{j+1}.}
\]
Since each triangle in the berry triangulation occurs in exactly three
necklaces of $R$, the equation 
\[
\sum_{N\text{ a necklace in }R}\text{Area}(N)=12\pi
\]
is satisfied.\vspace{2mm}
\item \textbf{The necklace combinatorial complementation constraint. }Let
$N$ be a necklace in the raspberry $R$. If $N$ has necklace code
$\necklace 0h{\ldots abc}\ldots$, then the raspberry $R$ must contain
a necklace that has necklace code $\necklace 0b{\ldots ahc}\ldots$.\vspace{2mm}
\end{enumerate}
\end{defn}

\section{Computation of necklace polynomials\label{sec:necklace-polynomials}}

Let $n\in\N$ and let $r_{1},\ldots,r_{n-1}$ be indeterminates and
$r_{0}:=1$. With $m\in\N$ and 
\[
\{h,b_{1},b_{2},\ldots,b_{m-1}\}\subseteq\{1,\ldots,n-1\},
\]
 consider the necklace code $\necklace 0h{b_{1}b_{2}\ldots b_{m-1}}$.
We can determine a polynomial equation in the indeterminates $r_{1},\ldots,r_{n-1}$
that the radii of spheres in a raspberry must necessarily satisfy
if the particular necklace code $\necklace 0h{b_{1}b_{2}\ldots b_{m-1}}$
appears as the code for a bona fide necklace in a raspberry. There
are various ways to compute this sought polynomial. None are particularly
fast. To compute these polynomials we employ Gr\"obner bases with
some precomputation. 

To this end, we introduce the following five blocks of symbols that
we regard as indeterminates for the purposes of the Gr\"obner basis
computations described in this section:
\begin{align*}
 & (\sin(\teddybear 0h{b_{j}b_{j+1}}))_{j\in\Z/m\Z}\\
 & (\cos(\teddybear 0h{b_{j}b_{j+1}}))_{j\in\Z/m\Z}\\
 & (\sin(\mickeymouse 0xy))_{x,y\in\{1,2,\ldots,n-1\}}\\
 & (\cos(\mickeymouse 0xy))_{x,y\in\{1,2,\ldots,n-1\}}\\
 & (r_{1},r_{2},\ldots,r_{n-1}).
\end{align*}
Each block of the above indeterminates is endowed with the degree
reverse lexicographic term order, and these orders are combined into
a block/product term order. This results in an elimination term order
so that a Gr\"obner basis computation of a set of polynomials with
respect to this term order favors eliminating indeterminates so that
only polynomials in $r_{1},r_{2},\ldots,r_{n-1}$ remain (cf. \cite[Section~2.3]{AdamsLoustaunau}
and \cite[Singular Manual, Section 3.3.3]{SINGULAR}).

Using the necklace dihedral angle constraint from Definition~\ref{def:constraints},
and applying the cosine on both sides, we obtain the equation 
\[
\cos\parenth{\sum_{j\in\Z/m\Z}\teddybear 0h{b_{j}b_{j+1}}}-1=0.
\]
Using angle addition identities, the left is expanded into a polynomial
equation in terms of sines and cosines of the expressions $\{\teddybear 0h{b_{j}b_{j+1}}\}_{j\in\Z/m\Z}$.
For the sake of brevity, we only write the unexpanded form. 

Consider the following set of polynomials which additionally capture
the relationship between indeterminates as given in Section~\ref{sec:Preliminaries}.
\begin{align*}
 & \left\{ \cos\parenth{\sum_{j\in\Z/m\Z}\teddybear 0h{b_{j}b_{j+1}}}-1\right\} \\
\cup & \left\{ \sin^{2}(\teddybear 0h{b_{j}b_{j+1}})+\cos^{2}(\teddybear 0h{b_{j}b_{j+1}})-1\right\} _{j\in\Z/m\Z}\\
\cup & \left\{ \begin{array}{l}
\cos(\teddybear 0h{b_{j}b_{j+1})}\sin(\mickeymouse 0h{b_{j}})\sin(\mickeymouse 0h{b_{j+1}})\\
\qquad-\cos(\mickeymouse 0{b_{j}}{b_{j+1}})+\cos(\mickeymouse 0h{b_{j}})\cos(\mickeymouse 0h{b_{j+1}})
\end{array}\right\} _{j\in\Z/m\Z}\\
\cup & \left\{ \sin^{2}(\mickeymouse 0xy)+\cos^{2}(\mickeymouse 0xy)-1\right\} _{x,y\in\{1,\ldots,n-1\}}\\
\cup & \left\{ \begin{array}{l}
2(r_{0}+r_{x})(r_{0}+r_{y})\cos(\mickeymouse 0xy)\\
\qquad-(r_{0}+r_{x})^{2}-(r_{0}+r_{y})^{2}+(r_{x}+r_{y})^{2}
\end{array}\right\} _{x,y\in\{1,\ldots,n-1\}}.
\end{align*}
Theoretically, computing a Gr\"obner basis of this set of polynomials
with the stated elimination order will yield a polynomial in only
$r_{1},\ldots,r_{n-1}$ in the computed Gr\"obner basis. This computation
is very time consuming and thus not practical. We get around this
difficulty by performing some precomputation and computing Gr\"obner
bases in two steps. 

As precomputation we eliminate all odd powers of sines from the expanded
form of the equation 
\[
\cos\parenth{\sum_{j\in\Z/m\Z}\teddybear 0h{b_{j}b_{j+1}}}-1=0.
\]
This is done by a grouping-and-squaring procedure (see for example
\cite[Section 5]{Fernique2021}). We give a general description: Group
terms of a multivariate polynomial equation onto two groups $X$ and
$Y$ (to be determined) satisfying $X+Y=0$. Replace the equation
$X+Y=0$ by the equation $X^{2}-Y^{2}=0$. Note that solution set
of the equation $X^{2}-Y^{2}=(X+Y)(X-Y)=0$ contains that of the original
equation $X+Y=0$. To determine the grouping, we see that, if $X$
contains all the terms with only odd powers of some fixed indeterminate
$u$, and $Y$ contains all the terms with only even powers of the
indeterminate $u$, then the new equation $X^{2}-Y^{2}=0$, when expanded,
is such that all terms contain only even powers of the indeterminate
$u$. For each indeterminate in a given set, this procedure is repeated,
until only even powers of the indeterminates in the set remain in
each term of the resulting polynomial. Applying this procedure, we
eliminate all odd powers of the indeterminates $\{\sin(\teddybear 0h{b_{j}b_{j+1}})\}_{j\in\Z/m\Z}$
from the above equation, which results in a polynomial equation $U=0$
in only cosines and even powers sines of the expressions $\{\teddybear 0h{b_{j}b_{j+1}}\}_{j\in\Z/m\Z}$. 

We compute a Gr\"obner basis with respect to the stated term order
for the set of polynomials below. Since no odd powers of the indeterminates
$\{\sin(\teddybear 0h{b_{j}b_{j+1}})\}_{j\in\Z/m\Z}$ occur in the
polynomial $U$, their elimination eased by the identity $\cos^{2}x+\sin^{2}x=1$.
\begin{align*}
 & \left\{ U\right\} \\
\cup & \left\{ \sin^{2}(\teddybear 0h{b_{j}b_{j+1}})+\cos^{2}(\teddybear 0h{b_{j}b_{j+1}})-1\right\} _{j\in\Z/m\Z}\\
\cup & \left\{ \begin{array}{l}
\cos(\teddybear 0h{b_{j}b_{j+1})}\sin(\mickeymouse 0h{b_{j}})\sin(\mickeymouse 0h{b_{j+1}})\\
\qquad-\cos(\mickeymouse 0{b_{j}}{b_{j+1}})+\cos(\mickeymouse 0h{b_{j}})\cos(\mickeymouse 0h{b_{j+1}})
\end{array}\right\} _{j\in\Z/m\Z}\\
\cup & \left\{ \sin^{2}(\mickeymouse 0xy)+\cos^{2}(\mickeymouse 0xy)-1\right\} _{x,y\in\{1,\ldots,n-1\}}.
\end{align*}
The Gr\"obner basis of the above set of polynomials contains a polynomial
$V$ in only the indeterminates $\left\{ \cos(\mickeymouse 0xy)\right\} _{x,y\in\{1,\ldots,n-1\}}.$ 

As a second step we compute a Gr\"obner basis with respect to the
stated term order for the set of polynomials:
\[
\left\{ V\right\} \cup\left\{ \begin{array}{l}
2(r_{0}+r_{x})(r_{0}+r_{y})\cos(\mickeymouse 0xy)\\
\qquad-(r_{0}+r_{x})^{2}-(r_{0}+r_{y})^{2}+(r_{x}+r_{y})^{2}
\end{array}\right\} _{x,y\in\{1,\ldots,n-1\}}.
\]
The resulting Gr\"obner basis contains a polynomial $W$ in only
the indeterminates $r_{1},r_{2},\ldots,r_{n-1}$.

The \emph{necklace polynomial} of the necklace code $\necklace 0h{b_{1}b_{2}\ldots b_{m-1}}$
is defined to be the irreducible factor of $W$ so that the solution
set, in $r_{1},\ldots,r_{n-1}$, of the equation 
\[
\sum_{j\in\Z/m\Z}\teddybear 0h{b_{j}b_{j+1}}=2\pi,
\]
is contained in the algebraic variety determined by this irreducible
factor of $W$.

\hspace{5mm}

We note that computing necklace polynomials is sometimes very time
and resource intensive. Indeed, for some 02-necklace codes computed
in Section~\ref{subsec:computing-pairs} we are unable to compute
the sought necklace polynomials using the procedure stated in this
section, as the computation (using Singular) exhausts 64GB of RAM
and is terminated before completion. This is luckily not a critical
failure, as explained in Section~\ref{subsec:excluding-necklace-code-pairs}.

We remark that an analogous result for necklaces to Descartes' Theorem
for flowers from \cite{MathewsZymaris} may improve the situation
by expediting the computation of necklace polynomials.

\section{Computing all raspberries with at most two sizes of berry\label{sec:Computing-all-raspberries}}

\subsection{The 01-necklace codes and necklace polynomials\label{subsec:The-01-necklace-codes-and-polynomials}}

Let $R$ be any raspberry with at most two sizes of berry with radii
$0<r_{1}\leq r_{2}$. We observe that for $\{x,y\}\subseteq\{1,2\}$
\[
\pi/3<\teddybear 01{xy}<\pi.
\]
The top bound of $\pi$ is strict because Definition~\ref{def:raspberry}(\ref{enu:no-hemispheres-or-lunes})
prohibits equality. The bottom bound $\pi/3$ is strict, as the infimum
$\pi/3$ over the domain $\{(r_{1},r_{2})\in\R^{2}:0<r_{1}\leq r_{2}\}$
is only taken on in the limit when $r_{2}=r_{1}$ in the limit $r_{2}\to0$. 

Given the necklace dihedral angle constraint from Definition~\ref{def:constraints},
all 01-necklace codes determined from a bona fide raspberry can thus
have only three, four or five beads. This allows us to enumerate all
such necklace codes and to compute the necklace polynomials for these
codes as described in Section~\ref{sec:necklace-polynomials}. 

If exactly two sizes of berry occur in a raspberry then two different
size berries must be tangent, and therefore a 01-necklace must determine
a 01-necklace code that contains a `$2$' as bead label (and vice
versa). We present all possible such 01-codes along with their necklace
polynomials in Table~\ref{tab:01-necklaces-and-polynomials}. 

\begin{table}[tbh]
\begin{tabular}{|c|l|l|}
\hline 
\noalign{\vskip1mm}
$\necklace 01{211}$ & \scalebox{0.7}{$-r_{1}r_{2}^{2}+r_{1}r_{2}-r_{1}+3r_{2}^{2}+3r_{2}$} & $r_{1}\in[3,\alpha]$\tabularnewline[1mm]
\hline 
\noalign{\vskip1mm}
$\necklace 01{221}$ & \scalebox{0.7}{$r_{1}^{2}r_{2}^{2}-2r_{1}^{2}r_{2}+r_{1}^{2}-2r_{1}r_{2}^{2}-4r_{1}r_{2}+r_{2}^{2}$} & $r_{1}\in[1,\alpha]$\tabularnewline[1mm]
\hline 
\noalign{\vskip1mm}
$\necklace 01{222}$ & \scalebox{0.7}{$-r_{1}^{2}r_{2}+3r_{1}^{2}+r_{1}r_{2}+3r_{1}-r_{2}$} & $r_{1}\in[0,\alpha]$\tabularnewline[1mm]
\hline 
\noalign{\vskip1mm}
$\necklace 01{2111}$ & \scalebox{0.7}{$-r_{1}^{2}r_{2}+r_{1}r_{2}+2r_{1}+r_{2}+1$} & $r_{1}\in[\phi,\beta]$\tabularnewline[1mm]
\hline 
\noalign{\vskip1mm}
$\necklace 01{2121}$ & \scalebox{0.7}{$-r_{1}r_{2}+r_{1}+r_{2}+1$} & $r_{1}\in[1,\beta]$\tabularnewline[1mm]
\hline 
\noalign{\vskip1mm}
$\necklace 01{2211}$ & \scalebox{0.7}{$r_{0}^{2}r_{1}^{2}-4r_{0}^{2}r_{1}r_{2}+r_{0}^{2}r_{2}^{2}-2r_{0}r_{1}^{2}r_{2}-2r_{0}r_{1}r_{2}^{2}+2r_{1}^{2}r_{2}^{2}$} & $r_{1}\in[0,\beta]$\tabularnewline[1mm]
\hline 
\noalign{\vskip1mm}
$\necklace 01{2221}$ & \scalebox{0.7}{$r_{1}^{3}r_{2}^{2}-2r_{1}^{3}r_{2}+r_{1}^{3}+r_{1}^{2}r_{2}^{2}-4r_{1}^{2}r_{2}+r_{1}^{2}-3r_{1}r_{2}+r_{2}^{2}$} & $r_{1}\in[0,\beta]$\tabularnewline[1mm]
\hline 
\noalign{\vskip1mm}
$\necklace 01{2222}$ & \scalebox{0.7}{$-r_{1}^{2}r_{2}+2r_{1}^{2}+2r_{1}-r_{2}$} & $r_{1}\in[0,\beta]$\tabularnewline[1mm]
\hline 
\noalign{\vskip1mm}
$\necklace 01{21111}$ & \scalebox{0.7}{$\begin{array}{l}
r_{1}^{4}r_{2}+r_{1}^{3}r_{2}^{2}-7r_{1}^{3}r_{2}-4r_{1}^{2}r_{2}^{2}-4r_{1}^{2}r_{2}+4r_{1}^{2}\\
\quad+4r_{1}r_{2}+4r_{1}+r_{2}^{2}+2r_{2}+1
\end{array}$} & $r_{1}\in[\psi,\gamma]$\tabularnewline[1mm]
\hline 
\noalign{\vskip1mm}
$\necklace 01{21211}$ & \scalebox{0.7}{$\begin{array}{l}
r_{1}^{4}r_{2}^{4}-4r_{1}^{4}r_{2}^{3}+6r_{1}^{4}r_{2}^{2}-4r_{1}^{4}r_{2}+r_{1}^{4}-6r_{1}^{3}r_{2}^{4}+8r_{1}^{3}r_{2}^{3}+6r_{1}^{3}r_{2}^{2}\\
\quad-8r_{1}^{3}r_{2}+2r_{1}^{3}+9r_{1}^{2}r_{2}^{4}+12r_{1}^{2}r_{2}^{3}-5r_{1}^{2}r_{2}^{2}-8r_{1}^{2}r_{2}+r_{1}^{2}-2r_{1}r_{2}^{4}\\
\quad-6r_{1}r_{2}^{3}-8r_{1}r_{2}^{2}-4r_{1}r_{2}+r_{2}^{4}+2r_{2}^{3}+r_{2}^{2}
\end{array}$} & $r_{1}\in[0,\gamma]$\tabularnewline[1mm]
\hline 
\noalign{\vskip1mm}
$\necklace 01{22111}$ & \scalebox{0.7}{$\begin{array}{l}
r_{1}^{8}r_{2}^{2}+2r_{1}^{7}r_{2}^{4}-6r_{1}^{7}r_{2}^{3}-4r_{1}^{7}r_{2}-14r_{1}^{6}r_{2}^{4}+18r_{1}^{6}r_{2}^{3}+2r_{1}^{6}r_{2}+4r_{1}^{6}\\
\quad+18r_{1}^{5}r_{2}^{4}+26r_{1}^{5}r_{2}^{3}-18r_{1}^{5}r_{2}^{2}-2r_{1}^{5}r_{2}+4r_{1}^{5}+25r_{1}^{4}r_{2}^{4}-20r_{1}^{4}r_{2}^{3}\\
\quad-21r_{1}^{4}r_{2}^{2}-6r_{1}^{4}r_{2}+r_{1}^{4}+14r_{1}^{3}r_{2}^{4}-34r_{1}^{3}r_{2}^{3}-6r_{1}^{3}r_{2}^{2}-2r_{1}^{3}r_{2}\\
\quad+12r_{1}^{2}r_{2}^{4}-14r_{1}^{2}r_{2}^{3}+6r_{1}r_{2}^{4}-2r_{1}r_{2}^{3}+r_{2}^{4}
\end{array}$} & $r_{1}\in[0,\gamma]$\tabularnewline[1mm]
\hline 
\noalign{\vskip1mm}
$\necklace 01{22121}$ & \scalebox{0.7}{$\begin{array}{l}
r_{1}^{5}r_{2}^{4}-4r_{1}^{5}r_{2}^{3}+6r_{1}^{5}r_{2}^{2}-4r_{1}^{5}r_{2}+r_{1}^{5}-6r_{1}^{4}r_{2}^{4}+8r_{1}^{4}r_{2}^{3}+4r_{1}^{4}r_{2}^{2}\\
\quad-8r_{1}^{4}r_{2}+2r_{1}^{4}+9r_{1}^{3}r_{2}^{4}+12r_{1}^{3}r_{2}^{3}-10r_{1}^{3}r_{2}^{2}-8r_{1}^{3}r_{2}+r_{1}^{3}\\
\quad-4r_{1}^{2}r_{2}^{3}-12r_{1}^{2}r_{2}^{2}-4r_{1}^{2}r_{2}+4r_{1}r_{2}^{4}+5r_{1}r_{2}^{3}+r_{2}^{4}+r_{2}^{3}
\end{array}$} & $r_{1}\in[0,\gamma]$\tabularnewline[1mm]
\hline 
\noalign{\vskip1mm}
$\necklace 01{22211}$ & \scalebox{0.7}{$\begin{array}{l}
r_{1}^{8}r_{2}^{5}-5r_{1}^{8}r_{2}^{4}+7r_{1}^{8}r_{2}^{3}-r_{1}^{8}r_{2}^{2}-2r_{1}^{8}r_{2}+r_{1}^{8}+r_{1}^{7}r_{2}^{6}-11r_{1}^{7}r_{2}^{5}\\
\quad+23r_{1}^{7}r_{2}^{4}-5r_{1}^{7}r_{2}^{3}+2r_{1}^{7}r_{2}^{2}-8r_{1}^{7}r_{2}+2r_{1}^{7}-3r_{1}^{6}r_{2}^{6}+63r_{1}^{6}r_{2}^{4}\\
\quad-64r_{1}^{6}r_{2}^{3}+13r_{1}^{6}r_{2}^{2}-12r_{1}^{6}r_{2}+r_{1}^{6}+6r_{1}^{5}r_{2}^{6}-6r_{1}^{5}r_{2}^{5}+60r_{1}^{5}r_{2}^{4}\\
\quad-104r_{1}^{5}r_{2}^{3}+22r_{1}^{5}r_{2}^{2}-6r_{1}^{5}r_{2}-4r_{1}^{4}r_{2}^{6}-21r_{1}^{4}r_{2}^{5}+64r_{1}^{4}r_{2}^{4}\\
\quad-73r_{1}^{4}r_{2}^{3}+13r_{1}^{4}r_{2}^{2}+3r_{1}^{3}r_{2}^{6}-7r_{1}^{3}r_{2}^{5}+57r_{1}^{3}r_{2}^{4}-25r_{1}^{3}r_{2}^{3}\\
\quad-12r_{1}^{2}r_{2}^{5}+22r_{1}^{2}r_{2}^{4}-8r_{1}r_{2}^{5}+r_{2}^{6}
\end{array}$} & $r_{1}\in[0,\gamma]$\tabularnewline[1mm]
\hline 
\noalign{\vskip1mm}
$\necklace 01{22221}$ & \scalebox{0.7}{$\begin{array}{l}
r_{1}^{6}r_{2}^{2}-2r_{1}^{6}r_{2}+r_{1}^{6}+2r_{1}^{5}r_{2}^{4}-10r_{1}^{5}r_{2}^{3}+12r_{1}^{5}r_{2}^{2}-6r_{1}^{5}r_{2}+2r_{1}^{5}\\
\quad-2r_{1}^{4}r_{2}^{4}-14r_{1}^{4}r_{2}^{3}+27r_{1}^{4}r_{2}^{2}-8r_{1}^{4}r_{2}+r_{1}^{4}-12r_{1}^{3}r_{2}^{3}+24r_{1}^{3}r_{2}^{2}\\
\quad-4r_{1}^{3}r_{2}-r_{1}^{2}r_{2}^{4}-12r_{1}^{2}r_{2}^{3}+9r_{1}^{2}r_{2}^{2}-6r_{1}r_{2}^{3}+r_{2}^{4}
\end{array}$} & $r_{1}\in[0,\gamma]$\tabularnewline[1mm]
\hline 
\noalign{\vskip1mm}
$\necklace 01{22222}$ & \scalebox{0.7}{$\begin{array}{l}
r_{1}^{4}r_{2}^{2}-5r_{1}^{4}r_{2}+5r_{1}^{4}-r_{1}^{3}r_{2}^{2}-5r_{1}^{3}r_{2}+10r_{1}^{3}+r_{1}^{2}r_{2}^{2}-5r_{1}^{2}r_{2}\\
\quad+5r_{1}^{2}-r_{1}r_{2}^{2}-5r_{1}r_{2}+r_{2}^{2}
\end{array}$} & $r_{1}\in[0,\gamma]$\tabularnewline[1mm]
\hline 
\end{tabular}\caption{\label{tab:01-necklaces-and-polynomials}All 01-necklace codes, associated
necklace polynomials and the allowed range for $r_{1}$ that can occur
in a raspberry with two sizes of berry. Here $\alpha\approx4.4494,$
$\beta\approx2.4142$, $\gamma\approx1.1085$ are roots of the polynomials
given in Figure~\ref{fig:raspberries-one-size} with $\phi:=(1+\sqrt{5})/2$
and $\psi\approx0.5374$ a root of $1-4x^{2}+x^{3}$.}
\end{table}

\subsection{The raspberries with one size of berry\label{subsec:raspberries-one-size}}

\begin{figure}[H]

\begin{tabular}{ccc}
\includegraphics[width=3cm]{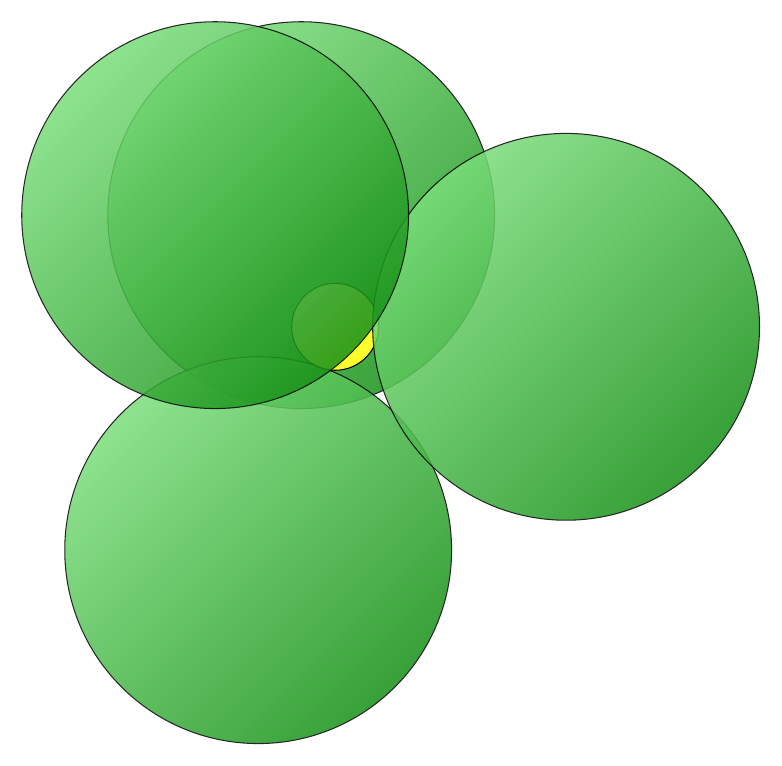} & \includegraphics[width=3cm]{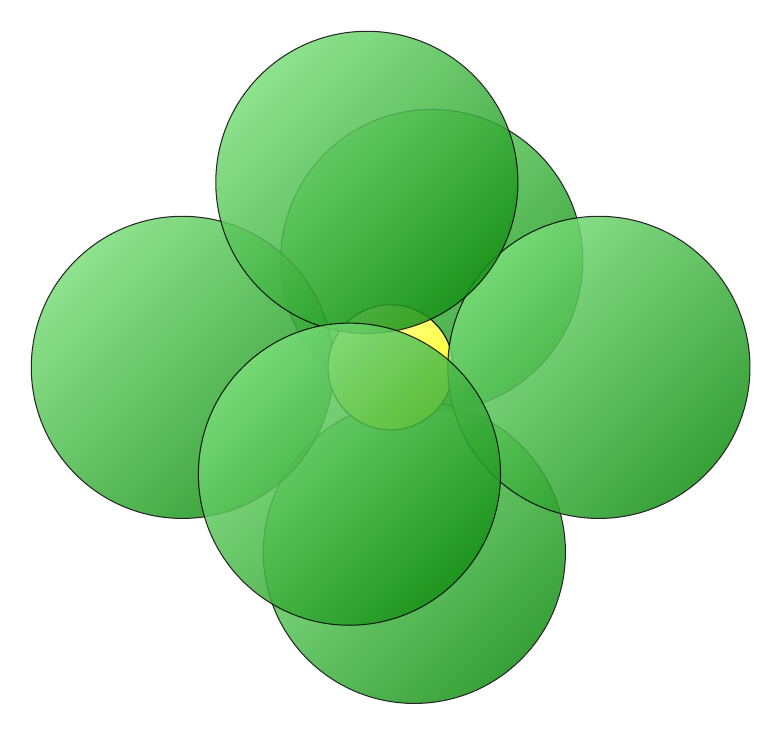} & \includegraphics[width=3cm]{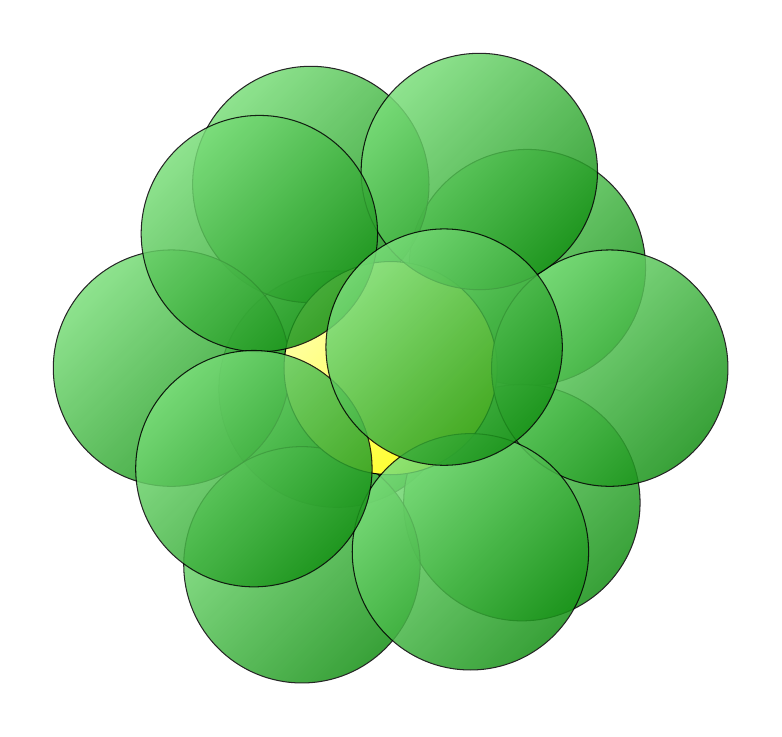}\tabularnewline
\scalebox{0.7}{$\begin{array}{l}
\necklace 01{111}\\
r_{1}\approx4.4494\\
r_{1}^{2}-4r_{1}-2
\end{array}$} & \scalebox{0.7}{$\begin{array}{l}
\necklace 01{111}1\\
r_{1}\approx2.4142\\
r_{1}^{2}-2r_{1}-1
\end{array}$} & \scalebox{0.7}{$\begin{array}{l}
\necklace 01{111}11\\
r_{1}\approx1.1085\\
r_{1}^{4}-6r_{1}^{3}+r_{1}^{2}+4r_{1}+1
\end{array}$}\tabularnewline
\end{tabular}\caption{The raspberries with one size of berry.\label{fig:raspberries-one-size}}
\end{figure}

All the raspberries with only one size of berry are displayed Figure~\ref{fig:raspberries-one-size}.
These are not surprising and have respective tetrahedral, octahedral,
and icosahedral symmetry and the radii of the berries are easily determined
from the necklace polynomials from the three possible 01-necklace
codes having only 1's as bead labels.

\subsection{The flexible raspberries with at most two sizes of berry\label{subsec:raspberries-flex}}

\begin{figure}[H]
\begin{tabular}{ccc}
\includegraphics[width=3cm]{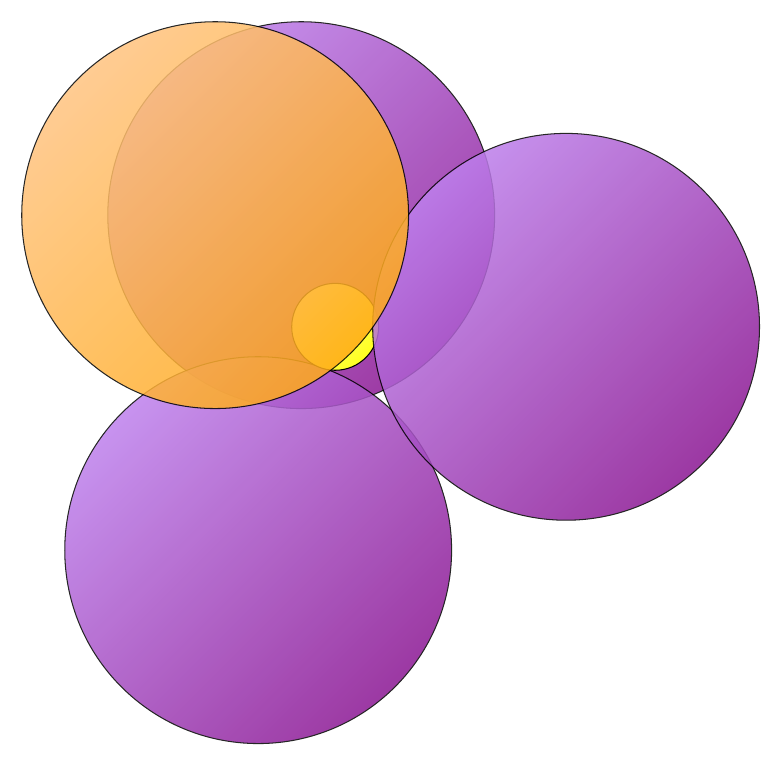} & \includegraphics[width=3cm]{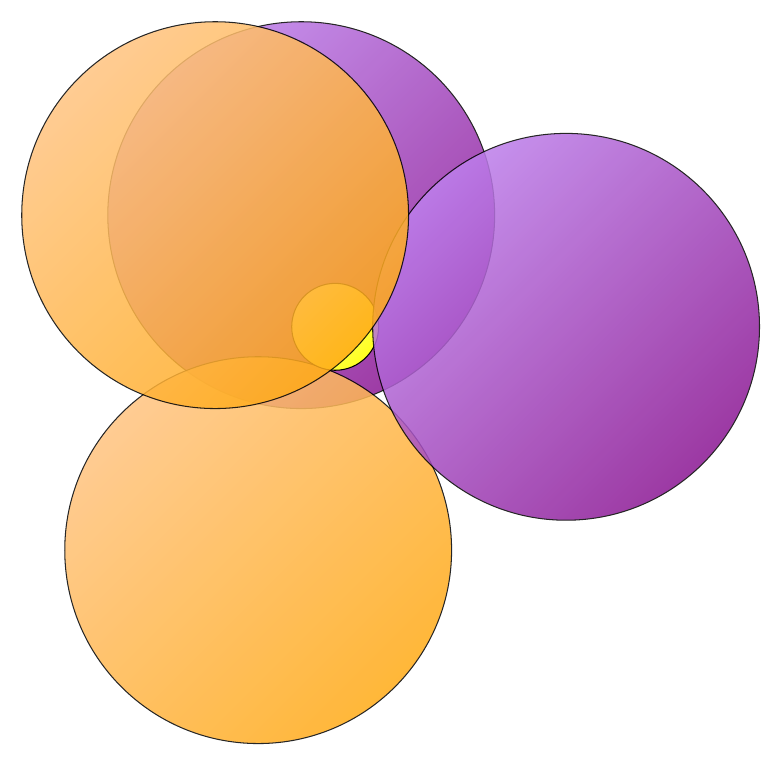} & \includegraphics[width=3cm]{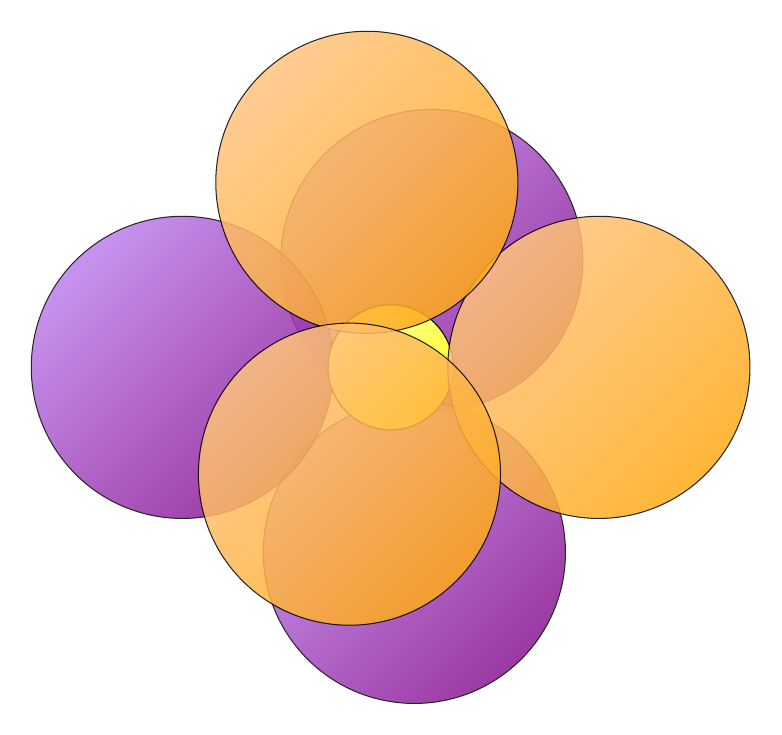}\tabularnewline
\scalebox{0.7}{$\begin{array}{l}
\necklace 01{222};\ \necklace 02{221}\\
-r_{1}^{2}r_{2}+3r_{1}^{2}+r_{1}r_{2}+3r_{1}-r_{2}
\end{array}$} & \scalebox{0.7}{$\begin{array}{l}
\necklace 01{221};\ \necklace 02{211}\\
r_{1}^{2}r_{2}^{2}-2r_{1}^{2}r_{2}+r_{1}^{2}\\
\quad-2r_{1}r_{2}^{2}-4r_{1}r_{2}+r_{2}^{2}
\end{array}$} & \scalebox{0.7}{$\begin{array}{l}
\necklace 01{2211};\ \necklace 02{2211}\\
2r_{1}^{2}r_{2}^{2}-2r_{1}^{2}r_{2}+r_{1}^{2}\\
\quad-2r_{1}r_{2}^{2}-4r_{1}r_{2}+r_{2}^{2}
\end{array}$}\tabularnewline
\end{tabular}\caption{The flexible raspberries with at most two sizes of berry.\label{fig:raspberries-flex}}
\end{figure}

All the raspberries that have at most two sizes of berry, but are
`flexible' are displayed in Figure~\ref{fig:raspberries-flex} along
with the necklace codes and necklace polynomials. The word `flexible'
here means that the berry radii is not uniquely determined by the
raspberry's combinatorics. This happens due to the 01-- and 02-necklace
codes arising from the raspberry, although different, both determine
the same necklace polynomial. 

That these are the only flexible raspberries that have at most two
sizes of berry is ex ante to the analysis done in Sections~\ref{subsec:Bounds-on-the-02-dihed-angles}--\ref{subsec:analyzing -numerical-computations}.
These are the only flexible raspberries that have at most two sizes
of berry that survive the culling process as described in Sections~\ref{subsec:Bounds-on-the-02-dihed-angles}--\ref{subsec:analyzing -numerical-computations}.

\subsection{Bounds on the dihedral angle $\protect\teddybear 02{xy}$ in a raspberry\label{subsec:Bounds-on-the-02-dihed-angles}}

The analysis presented in this section establishes that there can
be at most finitely many 02-necklace codes determined by 02-necklaces
in any raspberry with two sizes of berry. A necklace with arbitrarily
many (small) beads is easily imagined, hence, to be able to enumerate
all 02-necklace codes that can occur in a raspberry, we must determine
a strictly positive lower bound on each dihedral angle that can be
formed in a 02-necklace. The necklace dihedral angle constraint (Definition~\ref{def:constraints})
then limits the number of beads that can be present in such a necklace.
We may refer the reader to \cite[Section~4]{Fernique2021} where a
completely analogous computation is performed in two dimensions. However,
with the increase in dimension the situation is more subtle as it
requires optimization of the values $\teddybear 02{xy}$ with $\{x,y\}\subseteq\{1,2\}$
on algebraic varieties. \hspace{5mm}

Let $R$ be any raspberry with two sizes of berry with radii $0<r_{1}<r_{2}$.
Since we have computed all the possible 01-necklace codes that can
arise from $R$ (see Table~\ref{tab:01-necklaces-and-polynomials}),
the presence of a particular 01-necklace code places the constraint
on $(r_{1},r_{2})$ in that this point must lie on the algebraic variety
determined by the 01-necklace code's necklace polynomial. For each
01-necklace code, we compute upper and lower bounds for the values
of the dihedral angles $\teddybear 02{xy}$ with $\{x,y\}\subseteq\{1,2\}$
under these constraints. For the purposes of subsequently enumerating
all potential solutions to Section~\ref{subsec:computing-pairs}'s
Equation~(\ref{eq:02-necklace-angle-counts-equation}), obtaining
a strictly positive lower bound for the values of $\teddybear 02{xy}$
on this variety is important so that we are assured that only finitely
many potential solutions to Equation~(\ref{eq:02-necklace-angle-counts-equation})
exist. 

Using the data in Table~\ref{tab:01-necklaces-and-polynomials},
by discretizing the given interval in Table~\ref{tab:01-necklaces-and-polynomials},
with standard root finding techniques for univariate polynomials,
we compute a numerical discretization of the relevant part of the
algebraic variety determined by the necklace polynomial of each 01-necklace
code. Through visualization of the values of each expression $\teddybear 02{xy}$
with $\{x,y\}\subseteq\{1,2\}$ on these discretizations, it can be
seen that all are minimized where $r_{1}$ lies at one of the the
endpoints of the intervals from Table~\ref{tab:01-necklaces-and-polynomials}.
Since the exact values of $r_{1}$, and thus $r_{2}$, are known at
the endpoints, the minima can be determined exactly. 

Not all maxima of the values for $\teddybear 02{xy}$ with $\{x,y\}\subseteq\{1,2\}$
are determined where $r_{1}$ lies at one of endpoints of the intervals
in Table~\ref{tab:01-necklaces-and-polynomials} and as such, not
all maxima could be determined exactly. Maxima that could not be determined
exactly were estimated numerically through successive refinement of
the intervals of the varieties' discretizations that contain the maximum
until adequate precision is reached (here, to within $10^{-4}$ of
the true bound). 

For the purposes of computations, even if bounds could be determined
exactly, we use numerical approximations in which we conservatively
overestimate the upper bounds and underestimate the lower bounds,
so as to not inadvertently lose potential solutions to Section~\ref{subsec:computing-pairs}'s
Equation~(\ref{eq:02-necklace-angle-counts-equation}). The bounds
used in subsequent computations are displayed in Table~\ref{tab:02-teddy-bounds}.

\begin{table}
\begin{tabular}{|cl|ll|}
\hline 
\noalign{\vskip1mm}
$\necklace 01{211}$ & \scalebox{0.7}{$\begin{array}{ccll}
\teddybear 02{22} & \in & [2.094, & 3.142]\\
\teddybear 02{21} & \in & [2.094, & 3.142]\\
\teddybear 02{11} & \in & [2.094, & 2.095]
\end{array}$} & $\necklace 01{21111}$ & \scalebox{0.7}{$\begin{array}{cclc}
\teddybear 02{22} & \in & [1.256, & 3.142]\\
\teddybear 02{21} & \in & [0, & 1.257]\\
\teddybear 02{11} & \in & [0.750, & 1.257]
\end{array}$}\tabularnewline[1mm]
\hline 
\noalign{\vskip1mm}
$\necklace 01{221}$ & \scalebox{0.7}{$\begin{array}{ccll}
\teddybear 02{22} & \in & [2.094, & 3.142]\\
\teddybear 02{21} & \in & [2.094, & 2.618]\\
\teddybear 02{11} & \in & [1.047, & 2.095]
\end{array}$} & $\necklace 01{21211}$ & \scalebox{0.7}{$\begin{array}{cclc}
\teddybear 02{22} & \in & [1.047, & 1.319^{*}]\\
\teddybear 02{21} & \in & [0.261, & 1.257]\\
\teddybear 02{11} & \in & [0.523, & 1.257]
\end{array}$}\tabularnewline[1mm]
\hline 
\noalign{\vskip1mm}
$\necklace 01{222}$ & \scalebox{0.7}{$\begin{array}{ccll}
\teddybear 02{22} & \in & [1.047, & 3.142]\\
\teddybear 02{21} & \in & [0.523, & 2.095]\\
\teddybear 02{11} & \in & [0.585, & 2.095]
\end{array}$} & $\necklace 01{22111}$ & \scalebox{0.7}{$\begin{array}{cclc}
\teddybear 02{22} & \in & [1.047, & 1.279^{*}]\\
\teddybear 02{21} & \in & [0.710, & 1.257]\\
\teddybear 02{11} & \in & [0.672, & 1.257]
\end{array}$}\tabularnewline[1mm]
\hline 
\noalign{\vskip1mm}
$\necklace 01{2111}$ & \scalebox{0.7}{$\begin{array}{cclc}
\teddybear 02{22} & \in & [1.570, & 3.142]\\
\teddybear 02{21} & \in & [1.570, & 3.142]\\
\teddybear 02{11} & \in & [1.378, & 1.571]
\end{array}$} & $\necklace 01{22121}$ & \scalebox{0.7}{$\begin{array}{cclc}
\teddybear 02{22} & \in & [1.047, & 1.257]\\
\teddybear 02{21} & \in & [0.820, & 1.257]\\
\teddybear 02{11} & \in & [0.750, & 1.257]
\end{array}$}\tabularnewline[1mm]
\hline 
\noalign{\vskip1mm}
$\necklace 01{2121}$ & \scalebox{0.7}{$\begin{array}{cclc}
\teddybear 02{22} & \in & [1.570, & 3.142]\\
\teddybear 02{21} & \in & [1.570, & 1.571]\\
\teddybear 02{11} & \in & [1.047, & 1.571]
\end{array}$} & $\necklace 01{22211}$ & \scalebox{0.7}{$\begin{array}{cclc}
\teddybear 02{22} & \in & [1.0472, & 1.257]\\
\teddybear 02{21} & \in & [0.8523, & 1.257]\\
\teddybear 02{11} & \in & [0.7792, & 1.257]
\end{array}$}\tabularnewline[1mm]
\hline 
\noalign{\vskip1mm}
$\necklace 01{2211}$ & \scalebox{0.7}{$\begin{array}{cclc}
\teddybear 02{22} & \in & [1.047, & 3.142]\\
\teddybear 02{21} & \in & [0.261, & 1.571]\\
\teddybear 02{11} & \in & [0.523, & 1.571]
\end{array}$} & $\necklace 01{22221}$ & \scalebox{0.7}{$\begin{array}{cclc}
\teddybear 02{22} & \in & [1.047, & 1.257]\\
\teddybear 02{21} & \in & [0.907, & 1.257]\\
\teddybear 02{11} & \in & [0.836, & 1.257]
\end{array}$}\tabularnewline[1mm]
\hline 
\noalign{\vskip1mm}
$\necklace 01{2221}$ & \scalebox{0.7}{$\begin{array}{cclc}
\teddybear 02{22} & \in & [1.047, & 1.666]\\
\teddybear 02{21} & \in & [0.628, & 1.571]\\
\teddybear 02{11} & \in & [0.628, & 1.571]
\end{array}$} & $\necklace 01{22222}$ & \scalebox{0.7}{$\begin{array}{cclc}
\teddybear 02{22} & \in & [1.047, & 1.257]\\
\teddybear 02{21} & \in & [0.942, & 1.257]\\
\teddybear 02{11} & \in & [0.878, & 1.257]
\end{array}$}\tabularnewline[1mm]
\hline 
\noalign{\vskip1mm}
$\necklace 01{2222}$ & \scalebox{0.7}{$\begin{array}{cclc}
\teddybear 02{22} & \in & [1.047, & 1.571]\\
\teddybear 02{21} & \in & [0.785, & 1.571]\\
\teddybear 02{11} & \in & [0.722, & 1.571]
\end{array}$} &  & \multicolumn{1}{l}{}\tabularnewline[1mm]
\cline{1-2}
\end{tabular}\caption{\label{tab:02-teddy-bounds}Bounds on the values dihedral angles $\protect\teddybear 02{xy}$
with $\{x,y\}\subseteq\{1,2\}$ that can occur in a raspberry with
a given 01-necklace code present. Asterisks indicate bounds that could
not be determined exactly.}

\end{table}

\begin{rem}
\label{rem:21111-is-ok}The reader will notice the that the lower
bound of $\teddybear 02{21}$ for $\necklace 01{21111}$ is zero.
As mentioned in this section, obtaining a strictly positive lower
bound is important so that all possible solutions to Section~\ref{subsec:computing-pairs}'s
Equation~(\ref{eq:02-necklace-angle-counts-equation}) can be enumerated.
Here, this is actually not an issue. If the necklace code $\necklace 01{21111}$
arises from a raspberry, and a `$2$' occurs as a bead label in a
02-necklace code of such a raspberry, then there must exist a 01-necklace
with two adjacent $2$'s as bead labels which is covered by the other
01-necklace code cases. Therefore we only need to consider 02-necklace
codes with $1$'s as bead labels, and we only need the lower bound
for $\teddybear 02{11}$ in this case.
\end{rem}

\subsection{Computing compatible necklace code pairs\label{subsec:computing-pairs}}

Let $R$ be any raspberry with two sizes of berry with radii $0<r_{1}<r_{2}$.
Keeping Remark~\ref{rem:21111-is-ok} in mind, assuming a particular
01-necklace is present in the raspberry $R$, employing the corresponding
bounds in Table~\ref{tab:02-teddy-bounds} in interval arithmetic,
we compute, by brute force, all triples of non-negative integers $a,b,c$
that can potentially to solve the equation 
\begin{equation}
a\ \teddybear 02{11}+b\ \teddybear 02{21}+c\ \teddybear 02{22}=2\pi.\label{eq:02-necklace-angle-counts-equation}
\end{equation}
We compute at least one 02-necklace code that contains at least one
`$1$' as a bead label, at most five $2$'s as bead labels, and that
matches the counts of dihedral angles present in a potential solution
to Equation~(\ref{eq:02-necklace-angle-counts-equation}). If no
such 02-necklace code exists, we discard the solution. We record the
pairs of 01- and 02-necklace codes for which a solution potentially
exists. 

The number of such pairs for each 01-necklace code is displayed in
Table~\ref{tab:number-of-solutions}, totaling 596 pairs \cite[Dataset00.json]{raspberrydata}. 

\begin{table}
\begin{tabular}{|cl|ll|}
\hline 
\noalign{\vskip1mm}
$\necklace 01{211}$ & 3 & $\necklace 01{21111}$ & 4\tabularnewline[1mm]
\hline 
\noalign{\vskip1mm}
$\necklace 01{221}$ & 7 & $\necklace 01{21211}$ & 121\tabularnewline[1mm]
\hline 
\noalign{\vskip1mm}
$\necklace 01{222}$ & 76 & $\necklace 01{22111}$ & 43\tabularnewline[1mm]
\hline 
\noalign{\vskip1mm}
$\necklace 01{2111}$ & 7 & $\necklace 01{22121}$ & 30\tabularnewline[1mm]
\hline 
\noalign{\vskip1mm}
$\necklace 01{2121}$ & 9 & $\necklace 01{22211}$ & 27\tabularnewline[1mm]
\hline 
\noalign{\vskip1mm}
$\necklace 01{2211}$ & 127 & $\necklace 01{22221}$ & 22\tabularnewline[1mm]
\hline 
\noalign{\vskip1mm}
$\necklace 01{2221}$ & 62 & $\necklace 01{22222}$ & 19\tabularnewline[1mm]
\hline 
\noalign{\vskip1mm}
$\necklace 01{2222}$ & 39 &  & \multicolumn{1}{l}{}\tabularnewline[1mm]
\cline{1-2}
\end{tabular}

\caption{\label{tab:number-of-solutions}The number of potential solutions
to Eq. (\ref{eq:02-necklace-angle-counts-equation}) that match at
least one 02-necklace code with the given 01-necklace code present.}

\end{table}

\subsection{Excluding necklace code pairs\label{subsec:excluding-necklace-code-pairs}}

Let $R$ be any raspberry with two sizes of berry with radii $0<r_{1}<r_{2}$.
The total of 596 pairs of 01-- and 02-necklace codes (See Table~\ref{tab:number-of-solutions})
provides an upper bound for the number of pairs of radii $r_{1}$
and $r_{2}$ that can occur in non-flexible raspberries. For every
pair of necklace codes, we attempt to numerically solve for $r_{1}$
and $r_{2}$ in the corresponding system of necklace dihedral angle
constraints (see Definition~\ref{def:constraints}) determined by
the necklace codes in the pair. We obtain 539 pairs \cite[Dataset00.json, Label "A"]{raspberrydata}
resulting in a numerical solution and 57 pairs \cite[Dataset00.json, Labels "B" and "C"]{raspberrydata}
failing to find a numerical solution. 

The 57 pairs of necklace codes, for which no numerical solution is
obtained, must be carefully excluded. The obvious way to do this by
computing the necklace polynomials for the 02-necklace codes in these
pairs. Unfortunately this is not always possible, as is explained
at the end of Section~\ref{sec:necklace-polynomials}. Luckily this
is not a critical failure, as these pairs can be excluded by other
means.

Of these 57 pairs, 45 pairs \cite[Dataset00.json, Label "B"]{raspberrydata}
are such that, to satisfy the necklace combinatorial complementation
constraint (see Definition~\ref{def:constraints}), a different 01-necklace
code than the one in the pair must arise from a necklace in a hypothetical
raspberry. By computing the roots for all pairs necklace polynomials
for 01-necklace codes, we exclude each of these 45 pairs of necklace
codes by verifying that they cannot satisfy the necklace dihedral
angle constraints (see Definition~\ref{def:constraints}) when a
different 01-necklace is present in a hypothetical raspberry. 

For the remaining 12 pairs \cite[Dataset00.json, Label "C"]{raspberrydata}
of necklace codes we compute the necklace polynomials of their present
01-- and 02-necklace codes, compute their common roots, and verify
that they they cannot satisfy the necklace dihedral angle constraints
(see Definition~\ref{def:constraints}).

This excludes all of the 57 pairs of necklace codes for which no numerical
solution was found.

\subsection{Analyzing the numerical solutions to systems of necklace dihedral
angle constraints\label{subsec:analyzing -numerical-computations}}

From the remaining 539 pairs \cite[Dataset00.json, Label "A"]{raspberrydata}
of necklace codes for which numerical solutions in $(r_{1},r_{2})$
to the system of necklace dihedral angle constraints (see Definition~\ref{def:constraints})
were determined, we identify 354 unique\footnote{We regard numerical approximations as identical if they do not differ
by more than $10^{-20}.$} numerical solutions in $(r_{1},r_{2})$ \cite[Dataset01.json]{raspberrydata}.

Using interval arithmetic, for each of the 354 unique numerical solutions
in $(r_{1},r_{2})$, we compute all possible 01- and 02-necklace codes
that can satisfy the necklace dihedral angle constraint (see Definition~\ref{def:constraints})
with these approximate values of $r_{1}$ and $r_{2}$, and we test
whether the necklace combinatorial complementation constraint (see
Definition~\ref{def:constraints}) can be satisfied by a subset of
these necklace codes. Of the 354 solutions, 64 pairs of numerical
values $r_{1}$ and $r_{2}$ pass the test \cite[Dataset01.json, Labels "B"--"F"]{raspberrydata},
while 290 fail \cite[Dataset01.json, Label "A"]{raspberrydata}.

We thus have 64 pairs of numerical values for $r_{1}$ and $r_{2}$
along with an associated set of all necklace codes that can satisfy
the necklace dihedral angle constraint and the necklace combinatorial
complementation constraint (see Definition~\ref{def:constraints})
for these values of $r_{1}$ and $r_{2}$. For each of the 64 pairs
of values of $r_{1}$ and $r_{2}$ and each of the associated necklace
codes we compute the total areas of the spherical triangles that would
be formed on the pit by the head and beads of a necklace matching
the necklace code. For a set of necklace codes to be determined from
a bona fide raspberry, the pit area constraint (see Definition~\ref{def:constraints})
dictates that a non-negative integer linear combination of these areas
should add up to $12\pi$. Using interval arithmetic, of the 64 pairs
of numerical values for $r_{1}$ and $r_{2}$, we determine 14 pairs
that fail to have a solution \cite[Dataset01.json, Label "B"]{raspberrydata},
and 50 that succeed \cite[Dataset01.json, Labels "C"--"F"]{raspberrydata}.

We thus have 50 pairs of numerical values for $r_{1}$ and $r_{2}$
along with an associated set of all necklace codes. Of these, 12 pairs
\cite[Dataset01.json, Label "C"]{raspberrydata} are solutions to
one of the three flexible raspberries given in Section~\ref{subsec:raspberries-flex},
and one pair \cite[Dataset01.json, Label "D"]{raspberrydata} corresponds
to the icosahedral raspberry in Section~\ref{subsec:raspberries-one-size}
with $r_{1}=r_{2}$. 

\begin{figure}
\includegraphics[width=4cm]{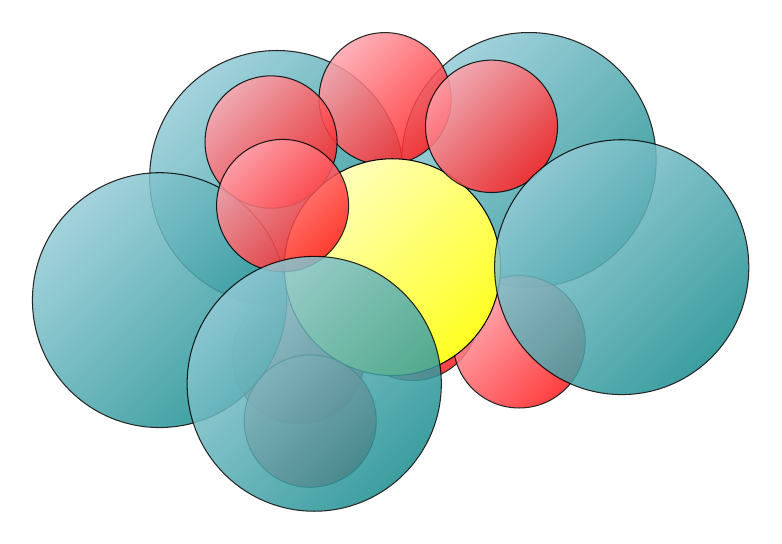}

\caption{The geometric failure in attempting to form a raspberry with necklace
codes $\protect\necklace 01{22111}$ and $\protect\necklace 02{211211}$
and $(r_{1},r_{2})\approx(0.6105,1.1754)$. The large berries must
occur tangent to an equator of the pit. However these berries are
too small for five of them to fit around the pit's equator, and too
large for six of them to fit around the pit's equator. \label{fig:geometric-failure}}
\end{figure}

Seven pairs of numerical values for $r_{1}$ and $r_{2}$ \cite[Dataset01.json, Label "E"]{raspberrydata}
fail to form raspberries. Of these seven, six pairs, presented in
Table~\ref{tab:tiling-failures}, are such that it is impossible
to construct a berry triangulation with vertices carrying labels from
$\{1,2\}$ so as to be compatible with the only possible 01- and 02-necklace
codes. This is easily verified by hand. Furthermore, of these seven,
one pair $(r_{1},r_{2})\approx(0.6105,1.1754)$ fails because its
geometry disallows forming a raspberry, see Figure~\ref{fig:geometric-failure}.
\begin{table}
\begin{tabular}{|c|c|}
\hline 
$(0.1248,0.5731)$ & $\necklace 01{21211};\necklace 02{11111111111}$\tabularnewline
\hline 
$(0.2869,1.1835)$ & $\necklace 01{21211};\necklace 02{111111111}$\tabularnewline
\hline 
$(0.3095,0.7962)$ & $\necklace 01{22111};\necklace 02{2211211}$\tabularnewline
\hline 
$(0.3263,1.1877)$ & $\necklace 01{2221};\necklace 02{21211211}$\tabularnewline
\hline 
$(0.3614,0.9015)$ & $\necklace 01{22111};\necklace 02{222211}$\tabularnewline
\hline 
$(0.4505,0.7999)$ & $\necklace 01{22121};\necklace 02{222111}$\tabularnewline
\hline 
\end{tabular}\caption{\label{tab:tiling-failures}The six numerical approximations for $(r_{1},r_{2})$
along with compatible necklace codes for which it is impossible to
construct a berry triangulation with vertices carrying labels from
$\{1,2\}$ so as to be compatible with the necklace codes.}
\end{table}

The remaining 30 pairs \cite[Dataset01.json, Label "F"]{raspberrydata}
of numerical values for $r_{1}$ and $r_{2}$ all determine non-flexible
raspberries. Examples of these we construct and exhibit in Section~\ref{sec:thirty-nonflexible}.

\section{The thirty non-flexible raspberries with two sizes of berry\label{sec:thirty-nonflexible}}

In this section we exhibit all thirty pairs of radii $(r_{1},r_{2})$
that permit the formation of a non-flexible raspberry with its berries
having radii $r_{1}$ and $r_{2}$. For each raspberry, we present
the necklace codes for all necklaces in the raspberry and approximations
of the values for $r_{1}$ and $r_{2}$, whose exact values are roots
of the given single variable polynomials. The raspberries are ordered
lexicographically by the necklace codes occurring in the raspberry
and indexed accordingly.

We note that, for given radii, the raspberries may not be unique.
For example, (24) permits the formation of two different raspberries
as discussed in \cite[Fig. 4]{FerniqueTwoSpheres2021} depending on
whether an antipodal pair of triangles formed by large berries are
aligned or misaligned through a $\pi/3$ rotation (the image below
is of the aligned version). Similarly, (25) also permits the formation
of two different raspberries where an antipodal pair of squares formed
by large berries are aligned or misaligned through a $\pi/4$ rotation
(the image below is of the misaligned version).

\noindent
\begin{minipage}[t]{0.60\textwidth} \vspace{0pt}    
    \noindent\textbf{
        \noindent (1) 
    }\noindent\necklace{0}{1}{21111};
\necklace{0}{2}{111111}

    \[
        (r_1, r_2) \approx (0.850, 1.940)
    \]

    \begin{dmath*}
    r_{1}^{4} - 8 r_{1}^{3} + 4 r_{1} + 1
    \end{dmath*}

    \begin{dmath*}
    2 r_{2}^{8} + 8 r_{2}^{7} + 4 r_{2}^{6} - 16 r_{2}^{5} - 40 r_{2}^{4} - 44 r_{2}^{3} - 26 r_{2}^{2} - 8 r_{2} - 1    
    \end{dmath*}

\end{minipage} \begin{minipage}[t]{0.35\textwidth} \vspace{0pt}
   \includegraphics[width=4cm]{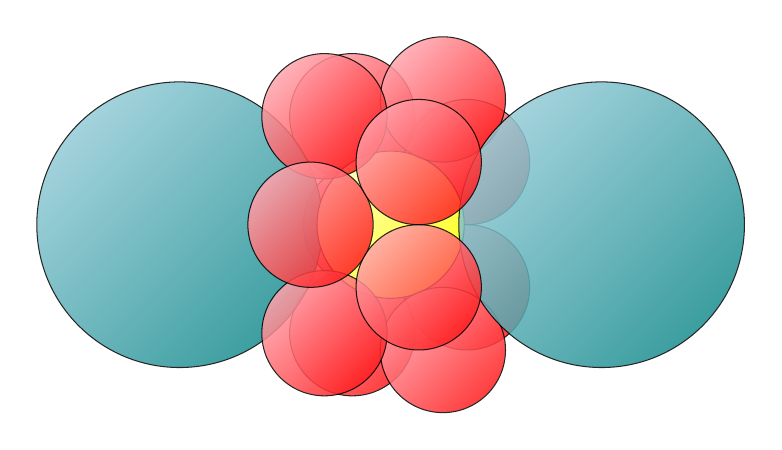}
\end{minipage}

\hspace{5mm}

\noindent
\begin{minipage}[t]{0.60\textwidth} \vspace{0pt}    
    \noindent\textbf{
        \noindent (2) 
    }\noindent\necklace{0}{1}{21111};
\necklace{0}{2}{1111111}

    \[
        (r_1, r_2) \approx (0.685, 4.012)
    \]

    \begin{dmath*}
    r_{1}^{6} - 10 r_{1}^{5} + 19 r_{1}^{4} + 16 r_{1}^{3} - 6 r_{1}^{2} - 6 r_{1} - 1
    \end{dmath*}

    \begin{dmath*}
    r_{2}^{12} + 40 r_{2}^{11} - 32 r_{2}^{10} - 452 r_{2}^{9} - 373 r_{2}^{8} - 300 r_{2}^{7} - 826 r_{2}^{6} - 888 r_{2}^{5} - 345 r_{2}^{4} + 10 r_{2}^{3} + 45 r_{2}^{2} + 12 r_{2} + 1    
    \end{dmath*}

\end{minipage} \begin{minipage}[t]{0.35\textwidth} \vspace{0pt}
   \includegraphics[width=4cm]{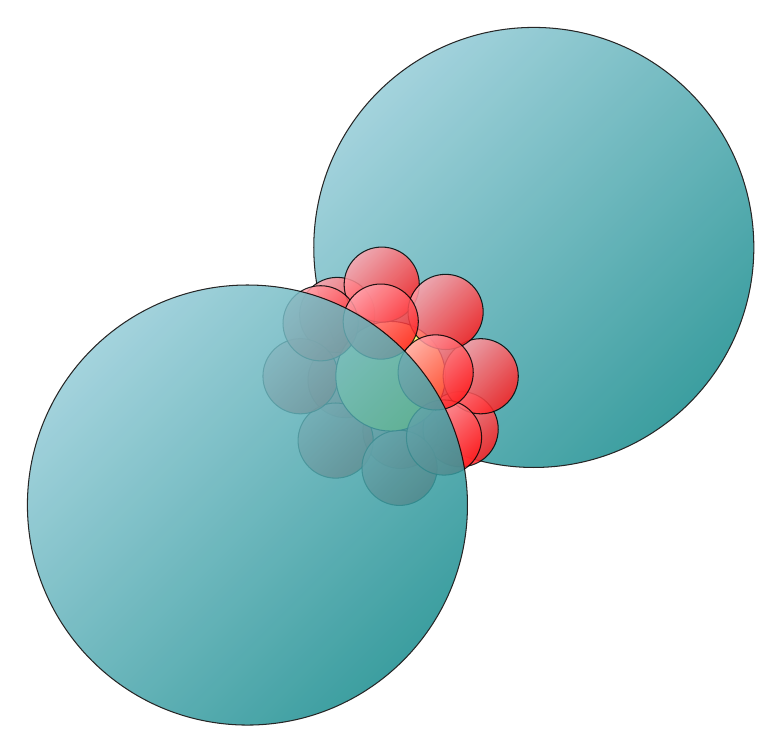}
\end{minipage}

\hspace{5mm}

\noindent
\begin{minipage}[t]{0.60\textwidth} \vspace{0pt}    
    \noindent\textbf{
        \noindent (3) 
    }\noindent\necklace{0}{1}{21111};
\necklace{0}{2}{11111111}

    \[
        (r_1, r_2) \approx (0.571, 17.31)
    \]

    \begin{dmath*}
    r_{1}^{8} - 16 r_{1}^{7} + 72 r_{1}^{6} - 48 r_{1}^{5} - 140 r_{1}^{4} - 32 r_{1}^{3} + 32 r_{1}^{2} + 16 r_{1} + 2
    \end{dmath*}

    \begin{dmath*}
    r_{2}^{16} - 272 r_{2}^{14} - 656 r_{2}^{13} + 2268 r_{2}^{12} + 11680 r_{2}^{11} + 21472 r_{2}^{10} + 20144 r_{2}^{9} + 7670 r_{2}^{8} - 3328 r_{2}^{7} - 5296 r_{2}^{6} - 2352 r_{2}^{5} - 196 r_{2}^{4} + 224 r_{2}^{3} + 96 r_{2}^{2} + 16 r_{2} + 1    
    \end{dmath*}

\end{minipage} \begin{minipage}[t]{0.35\textwidth} \vspace{0pt}
   \includegraphics[width=4cm]{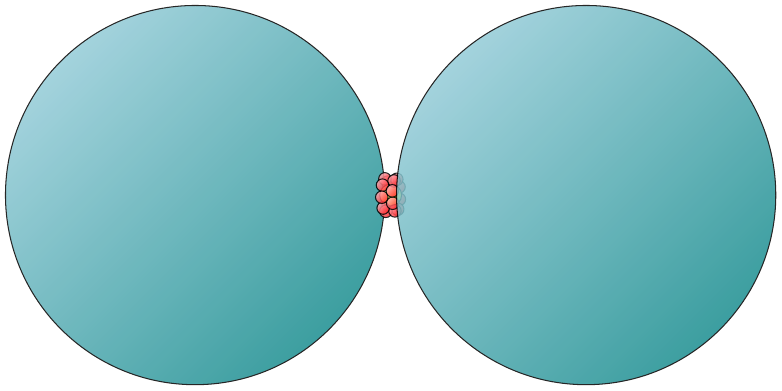}
\end{minipage}

\hspace{5mm}

\noindent
\begin{minipage}[t]{0.60\textwidth} \vspace{0pt}    
    \noindent\textbf{
        \noindent (4) 
    }\noindent\necklace{0}{1}{21111};
\necklace{0}{1}{2211};
\necklace{0}{2}{22111}

    \[
        (r_1, r_2) \approx (0.633, 6.107)
    \]

    \begin{dmath*}
    2 r_{1}^{4} - 12 r_{1}^{3} - 2 r_{1}^{2} + 4 r_{1} + 1
    \end{dmath*}

    \begin{dmath*}
    2 r_{2}^{4} - 12 r_{2}^{3} - 2 r_{2}^{2} + 4 r_{2} + 1    
    \end{dmath*}

\end{minipage} \begin{minipage}[t]{0.35\textwidth} \vspace{0pt}
   \includegraphics[width=4cm]{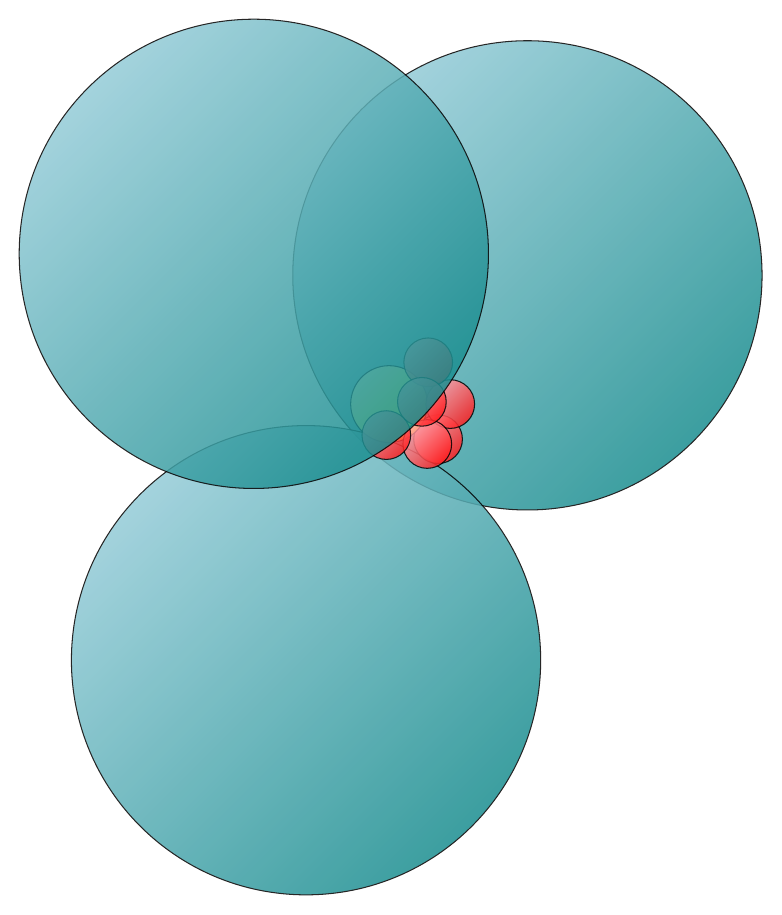}
\end{minipage}

\hspace{5mm}

\noindent
\begin{minipage}[t]{0.60\textwidth} \vspace{0pt}    
    \noindent\textbf{
        \noindent (5) 
    }\noindent\necklace{0}{1}{2121};
\necklace{0}{2}{11111}

    \[
        (r_1, r_2) \approx (1.425, 5.695)
    \]

    \begin{dmath*}
    r_{1}^{4} - 20 r_{1}^{3} + 10 r_{1}^{2} + 20 r_{1} + 5
    \end{dmath*}

    \begin{dmath*}
    r_{2}^{4} - 6 r_{2}^{3} + r_{2}^{2} + 4 r_{2} + 1    
    \end{dmath*}

\end{minipage} \begin{minipage}[t]{0.35\textwidth} \vspace{0pt}
   \includegraphics[width=4cm]{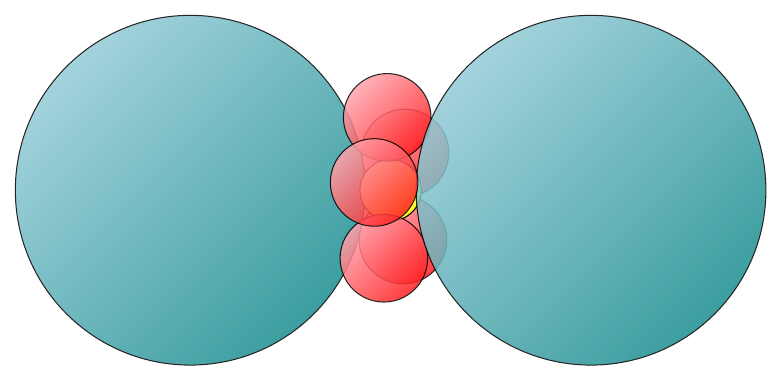}
\end{minipage}

\hspace{5mm}

\noindent
\begin{minipage}[t]{0.60\textwidth} \vspace{0pt}    
    \noindent\textbf{
        \noindent (6) 
    }\noindent\necklace{0}{1}{21211};
\necklace{0}{2}{111111}

    \[
        (r_1, r_2) \approx (0.743, 1.273)
    \]

    \begin{dmath*}
    9 r_{1}^{2} - 4 r_{1} - 2
    \end{dmath*}

    \begin{dmath*}
    3 r_{2}^{4} + 6 r_{2}^{3} - 5 r_{2}^{2} - 8 r_{2} - 2    
    \end{dmath*}

\end{minipage} \begin{minipage}[t]{0.35\textwidth} \vspace{0pt}
   \includegraphics[width=4cm]{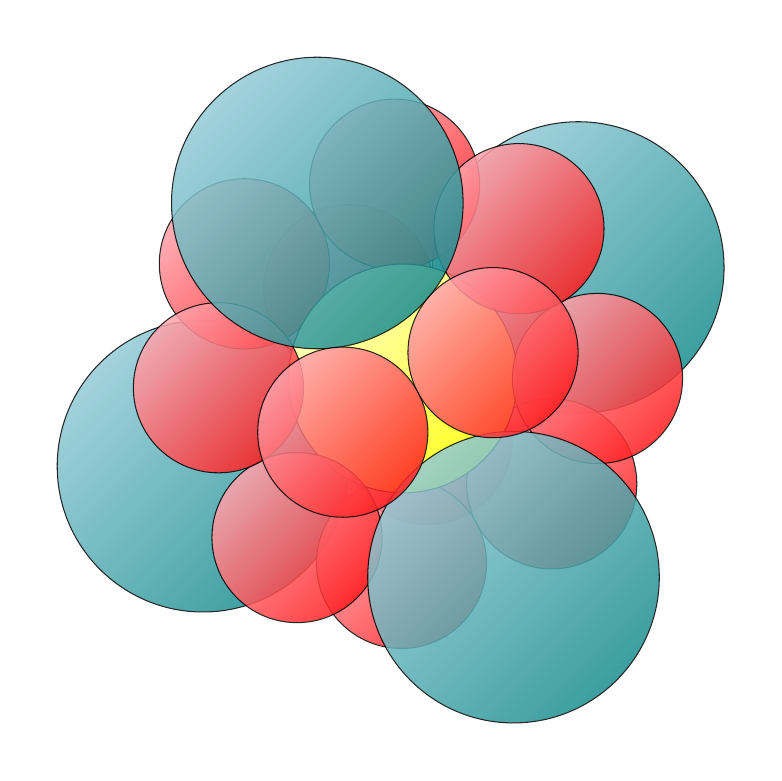}
\end{minipage}

\hspace{5mm}

\noindent
\begin{minipage}[t]{0.60\textwidth} \vspace{0pt}    
    \noindent\textbf{
        \noindent (7) 
    }\noindent\necklace{0}{1}{21211};
\necklace{0}{2}{11111111}

    \[
        (r_1, r_2) \approx (0.390, 1.334)
    \]

    \begin{dmath*}
    4 r_{1}^{4} - 24 r_{1}^{3} - 8 r_{1}^{2} + 4 r_{1} + 1
    \end{dmath*}

    \begin{dmath*}
    r_{2}^{4} - 4 r_{2}^{3} + 4 r_{2} + 1    
    \end{dmath*}

\end{minipage} \begin{minipage}[t]{0.35\textwidth} \vspace{0pt}
   \includegraphics[width=4cm]{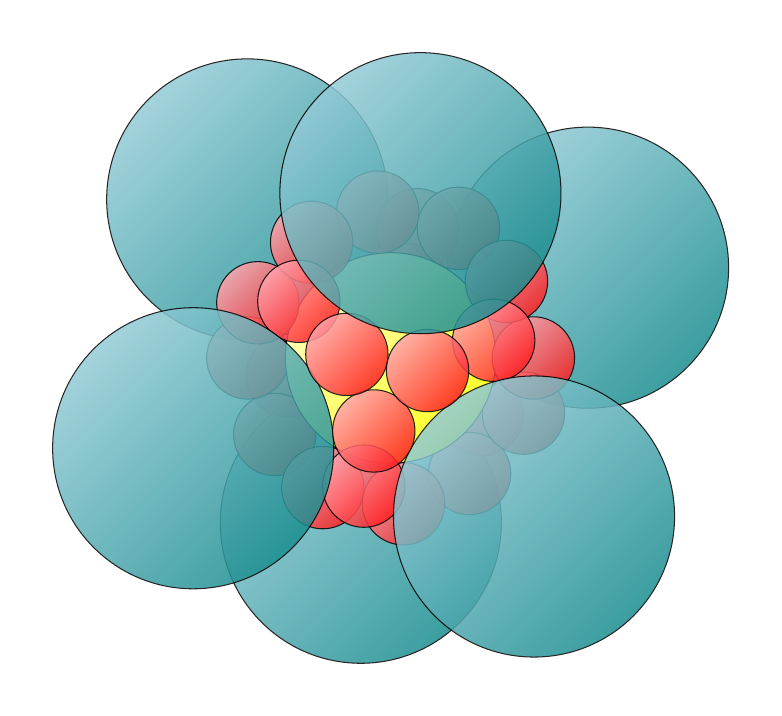}
\end{minipage}

\hspace{5mm}

\noindent
\begin{minipage}[t]{0.60\textwidth} \vspace{0pt}    
    \noindent\textbf{
        \noindent (8) 
    }\noindent\necklace{0}{1}{21211};
\necklace{0}{2}{1111111111}

    \[
        (r_1, r_2) \approx (0.202, 0.921)
    \]

    \begin{dmath*}
    25 r_{1}^{4} - 70 r_{1}^{3} - 31 r_{1}^{2} + 4 r_{1} + 1
    \end{dmath*}

    \begin{dmath*}
    4 r_{2}^{8} - 24 r_{2}^{7} - 48 r_{2}^{6} + 20 r_{2}^{5} + 59 r_{2}^{4} + 10 r_{2}^{3} - 17 r_{2}^{2} - 8 r_{2} - 1    
    \end{dmath*}

\end{minipage} \begin{minipage}[t]{0.35\textwidth} \vspace{0pt}
   \includegraphics[width=4cm]{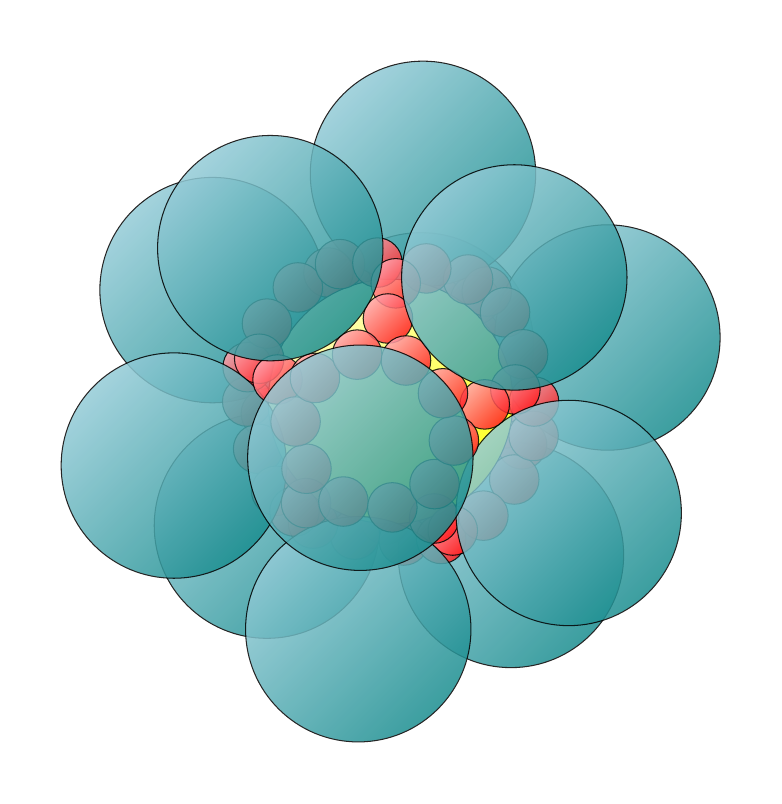}
\end{minipage}

\hspace{5mm}

\noindent
\begin{minipage}[t]{0.60\textwidth} \vspace{0pt}    
    \noindent\textbf{
        \noindent (9) 
    }\noindent\necklace{0}{1}{2211};
\necklace{0}{2}{211211}

    \[
        (r_1, r_2) \approx (0.763, 6.464)
    \]

    \begin{dmath*}
    13 r_{1}^{2} - 6 r_{1} - 3
    \end{dmath*}

    \begin{dmath*}
    r_{2}^{2} - 6 r_{2} - 3    
    \end{dmath*}

\end{minipage} \begin{minipage}[t]{0.35\textwidth} \vspace{0pt}
   \includegraphics[width=4cm]{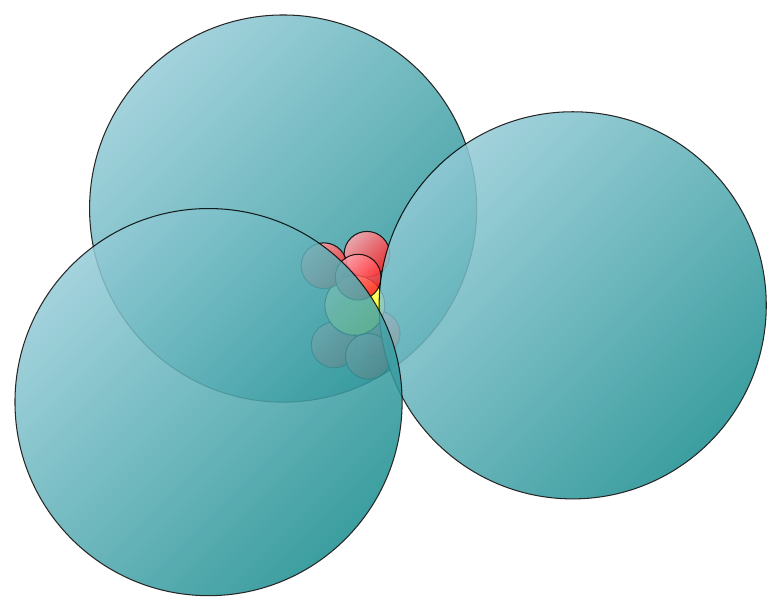}
\end{minipage}

\hspace{5mm}

\noindent
\begin{minipage}[t]{0.60\textwidth} \vspace{0pt}    
    \noindent\textbf{
        \noindent (10) 
    }\noindent\necklace{0}{1}{2211};
\necklace{0}{2}{211211211}

    \[
        (r_1, r_2) \approx (0.463, 4.449)
    \]

    \begin{dmath*}
    73 r_{1}^{4} - 88 r_{1}^{3} - 28 r_{1}^{2} + 16 r_{1} + 4
    \end{dmath*}

    \begin{dmath*}
    r_{2}^{2} - 4 r_{2} - 2    
    \end{dmath*}

\end{minipage} \begin{minipage}[t]{0.35\textwidth} \vspace{0pt}
   \includegraphics[width=4cm]{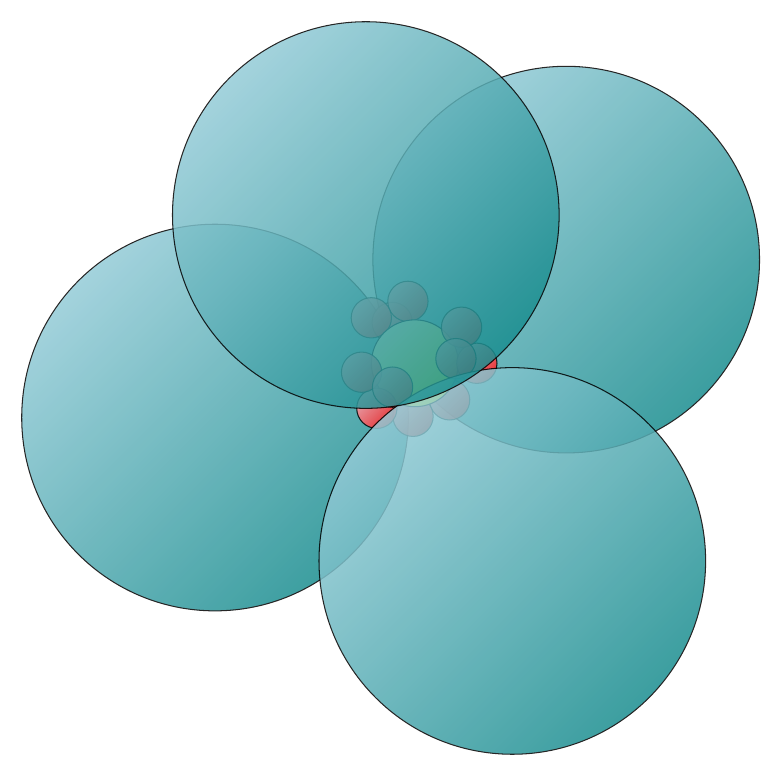}
\end{minipage}

\hspace{5mm}

\noindent
\begin{minipage}[t]{0.60\textwidth} \vspace{0pt}    
    \noindent\textbf{
        \noindent (11) 
    }\noindent\necklace{0}{1}{2211};
\necklace{0}{2}{211211211211}

    \[
        (r_1, r_2) \approx (0.308, 2.414)
    \]

    \begin{dmath*}
    17 r_{1}^{2} - 2 r_{1} - 1
    \end{dmath*}

    \begin{dmath*}
    r_{2}^{2} - 2 r_{2} - 1    
    \end{dmath*}

\end{minipage} \begin{minipage}[t]{0.35\textwidth} \vspace{0pt}
   \includegraphics[width=4cm]{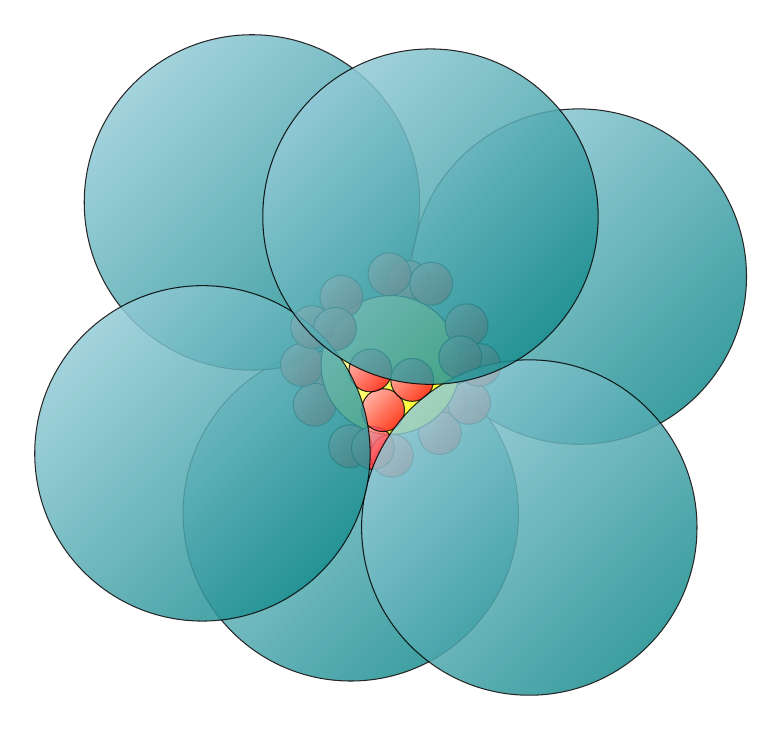}
\end{minipage}

\hspace{5mm}

\noindent
\begin{minipage}[t]{0.60\textwidth} \vspace{0pt}    
    \noindent\textbf{
        \noindent (12) 
    }\noindent\necklace{0}{1}{2211};
\necklace{0}{2}{211211211211211}

    \[
        (r_1, r_2) \approx (0.184, 1.108)
    \]

    \begin{dmath*}
    61 r_{1}^{4} + 2 r_{1}^{3} - 21 r_{1}^{2} - 2 r_{1} + 1
    \end{dmath*}

    \begin{dmath*}
    r_{2}^{4} - 6 r_{2}^{3} + r_{2}^{2} + 4 r_{2} + 1    
    \end{dmath*}

\end{minipage} \begin{minipage}[t]{0.35\textwidth} \vspace{0pt}
   \includegraphics[width=4cm]{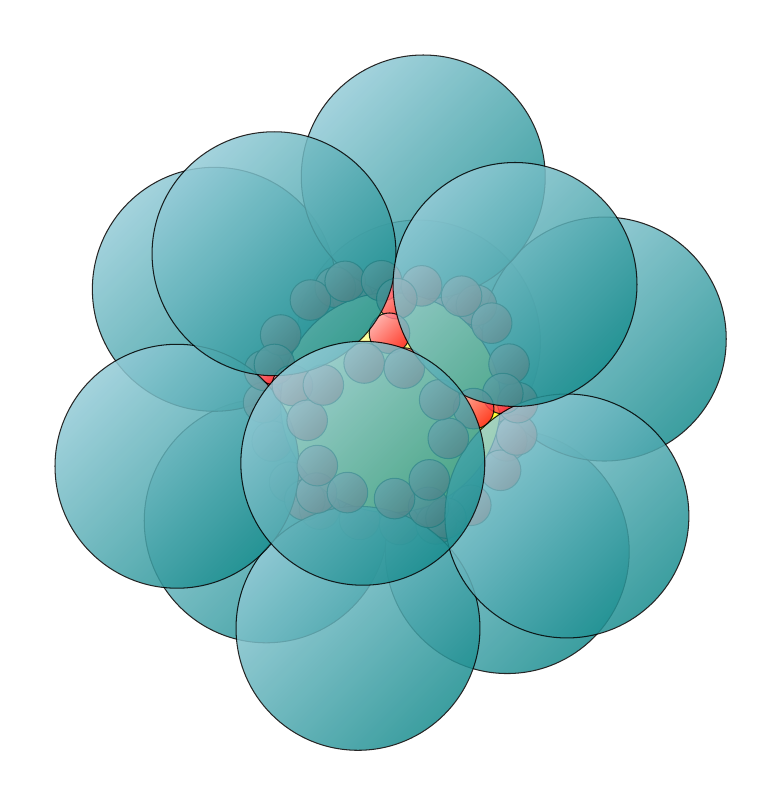}
\end{minipage}

\hspace{5mm}

\noindent
\begin{minipage}[t]{0.60\textwidth} \vspace{0pt}    
    \noindent\textbf{
        \noindent (13) 
    }\noindent\necklace{0}{1}{2211};
\necklace{0}{1}{222};
\necklace{0}{2}{21211}

    \[
        (r_1, r_2) \approx (0.976, 5.925)
    \]

    \begin{dmath*}
    13 r_{1}^{4} + 4 r_{1}^{3} - 6 r_{1}^{2} - 8 r_{1} - 2
    \end{dmath*}

    \begin{dmath*}
    r_{2}^{4} - 4 r_{2}^{3} - 10 r_{2}^{2} - 8 r_{2} - 2    
    \end{dmath*}

\end{minipage} \begin{minipage}[t]{0.35\textwidth} \vspace{0pt}
   \includegraphics[width=4cm]{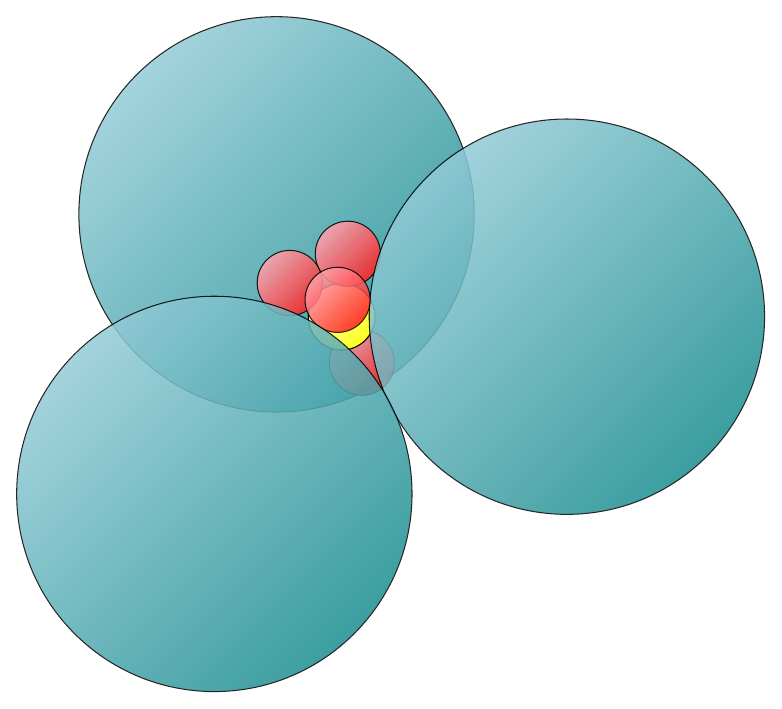}
\end{minipage}

\hspace{5mm}

\noindent
\begin{minipage}[t]{0.60\textwidth} \vspace{0pt}    
    \noindent\textbf{
        \noindent (14) 
    }\noindent\necklace{0}{1}{222};
\necklace{0}{2}{2121}

    \[
        (r_1, r_2) \approx (1.366, 6.464)
    \]

    \begin{dmath*}
    2 r_{1}^{2} - 2 r_{1} - 1
    \end{dmath*}

    \begin{dmath*}
    r_{2}^{2} - 6 r_{2} - 3    
    \end{dmath*}

\end{minipage} \begin{minipage}[t]{0.35\textwidth} \vspace{0pt}
   \includegraphics[width=4cm]{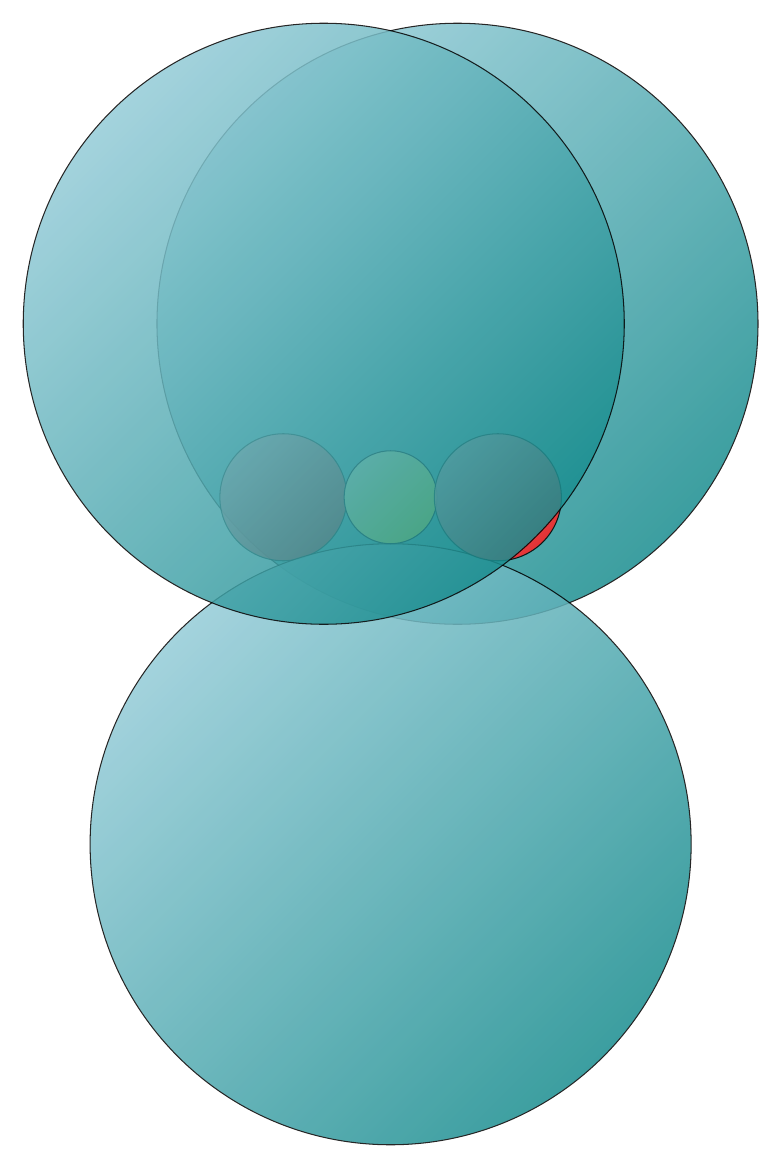}
\end{minipage}

\hspace{5mm}

\noindent
\begin{minipage}[t]{0.60\textwidth} \vspace{0pt}    
    \noindent\textbf{
        \noindent (15) 
    }\noindent\necklace{0}{1}{222};
\necklace{0}{2}{212121}

    \[
        (r_1, r_2) \approx (0.689, 4.449)
    \]

    \begin{dmath*}
    5 r_{1}^{2} - 2 r_{1} - 1
    \end{dmath*}

    \begin{dmath*}
    r_{2}^{2} - 4 r_{2} - 2    
    \end{dmath*}

\end{minipage} \begin{minipage}[t]{0.35\textwidth} \vspace{0pt}
   \includegraphics[width=4cm]{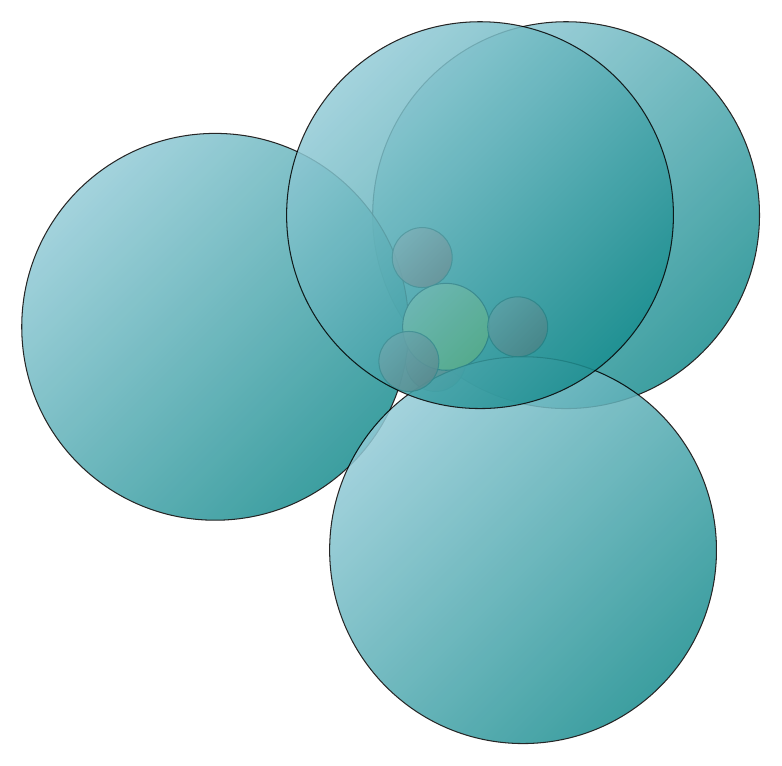}
\end{minipage}

\hspace{5mm}

\noindent
\begin{minipage}[t]{0.60\textwidth} \vspace{0pt}    
    \noindent\textbf{
        \noindent (16) 
    }\noindent\necklace{0}{1}{222};
\necklace{0}{2}{21212121}

    \[
        (r_1, r_2) \approx (0.426, 2.414)
    \]

    \begin{dmath*}
    2 r_{1}^{4} + 20 r_{1}^{3} + 6 r_{1}^{2} - 4 r_{1} - 1
    \end{dmath*}

    \begin{dmath*}
    r_{2}^{2} - 2 r_{2} - 1    
    \end{dmath*}

\end{minipage} \begin{minipage}[t]{0.35\textwidth} \vspace{0pt}
   \includegraphics[width=4cm]{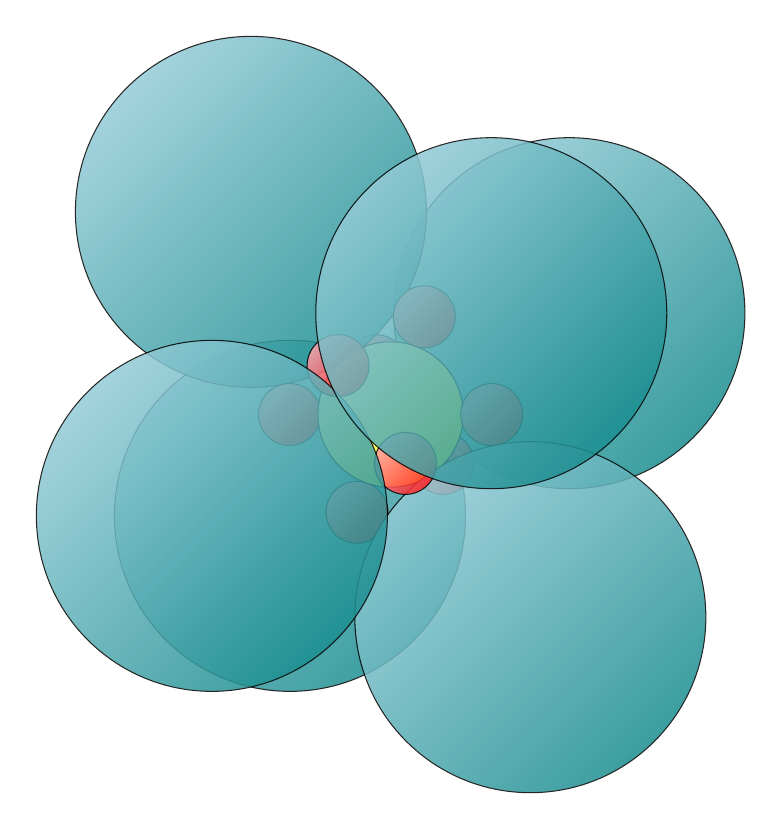}
\end{minipage}

\hspace{5mm}

\noindent
\begin{minipage}[t]{0.60\textwidth} \vspace{0pt}    
    \noindent\textbf{
        \noindent (17) 
    }\noindent\necklace{0}{1}{222};
\necklace{0}{2}{2121212121}

    \[
        (r_1, r_2) \approx (0.242, 1.108)
    \]

    \begin{dmath*}
    59 r_{1}^{8} + 28 r_{1}^{7} - 370 r_{1}^{6} - 152 r_{1}^{5} + 236 r_{1}^{4} + 142 r_{1}^{3} + 5 r_{1}^{2} - 8 r_{1} - 1
    \end{dmath*}

    \begin{dmath*}
    r_{2}^{4} - 6 r_{2}^{3} + r_{2}^{2} + 4 r_{2} + 1    
    \end{dmath*}

\end{minipage} \begin{minipage}[t]{0.35\textwidth} \vspace{0pt}
   \includegraphics[width=4cm]{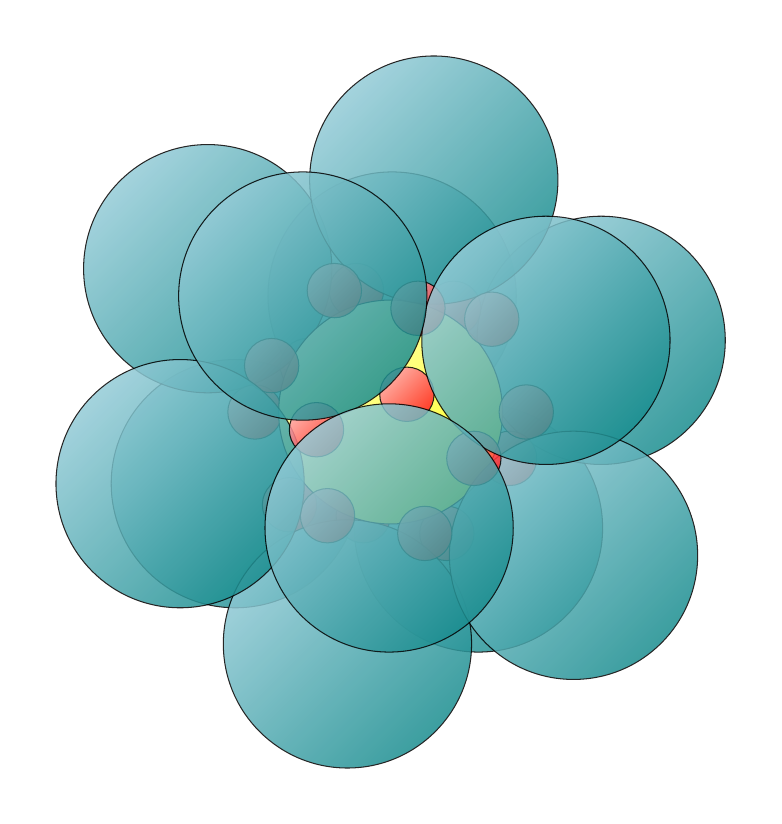}
\end{minipage}

\hspace{5mm}

\noindent
\begin{minipage}[t]{0.60\textwidth} \vspace{0pt}    
    \noindent\textbf{
        \noindent (18) 
    }\noindent\necklace{0}{1}{2221};
\necklace{0}{2}{21211}

    \[
        (r_1, r_2) \approx (1.052, 2.780)
    \]

    \begin{dmath*}
    9 r_{1}^{4} + 16 r_{1}^{3} - 8 r_{1}^{2} - 16 r_{1} - 4
    \end{dmath*}

    \begin{dmath*}
    9 r_{2}^{4} - 28 r_{2}^{3} + 2 r_{2}^{2} + 16 r_{2} + 4    
    \end{dmath*}

\end{minipage} \begin{minipage}[t]{0.35\textwidth} \vspace{0pt}
   \includegraphics[width=4cm]{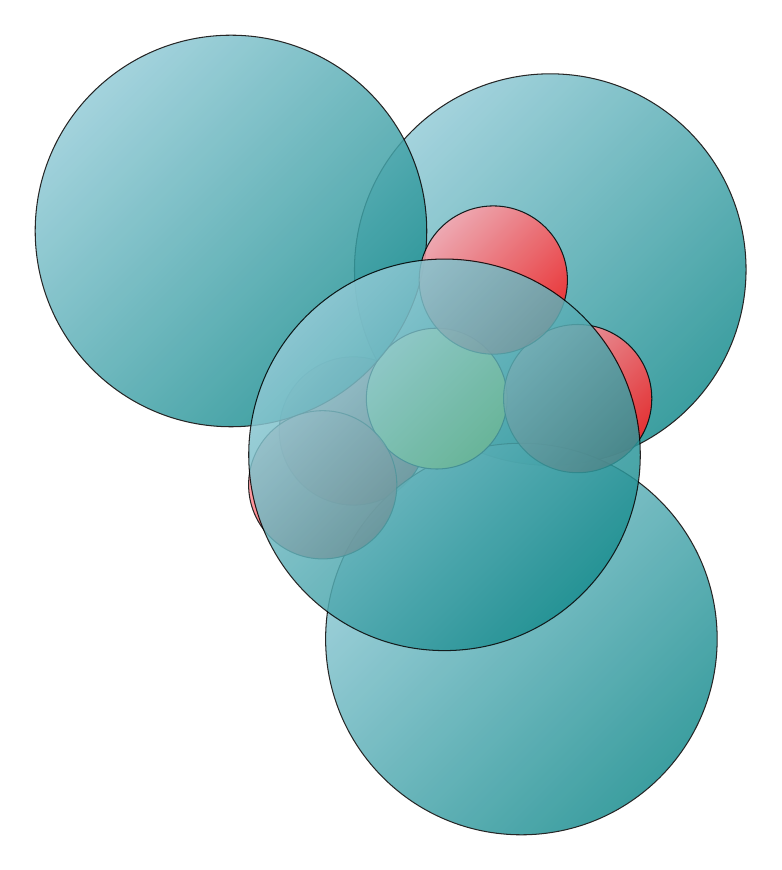}
\end{minipage}

\hspace{5mm}

\noindent
\begin{minipage}[t]{0.60\textwidth} \vspace{0pt}    
    \noindent\textbf{
        \noindent (19) 
    }\noindent\necklace{0}{1}{2221};
\necklace{0}{2}{221211}

    \[
        (r_1, r_2) \approx (0.554, 2.028)
    \]

    \begin{dmath*}
    36 r_{1}^{6} + 178 r_{1}^{5} + 281 r_{1}^{4} + 120 r_{1}^{3} - 60 r_{1}^{2} - 54 r_{1} - 9
    \end{dmath*}

    \begin{dmath*}
    81 r_{2}^{6} - 230 r_{2}^{5} + 89 r_{2}^{4} + 132 r_{2}^{3} - 57 r_{2}^{2} - 54 r_{2} - 9    
    \end{dmath*}

\end{minipage} \begin{minipage}[t]{0.35\textwidth} \vspace{0pt}
   \includegraphics[width=4cm]{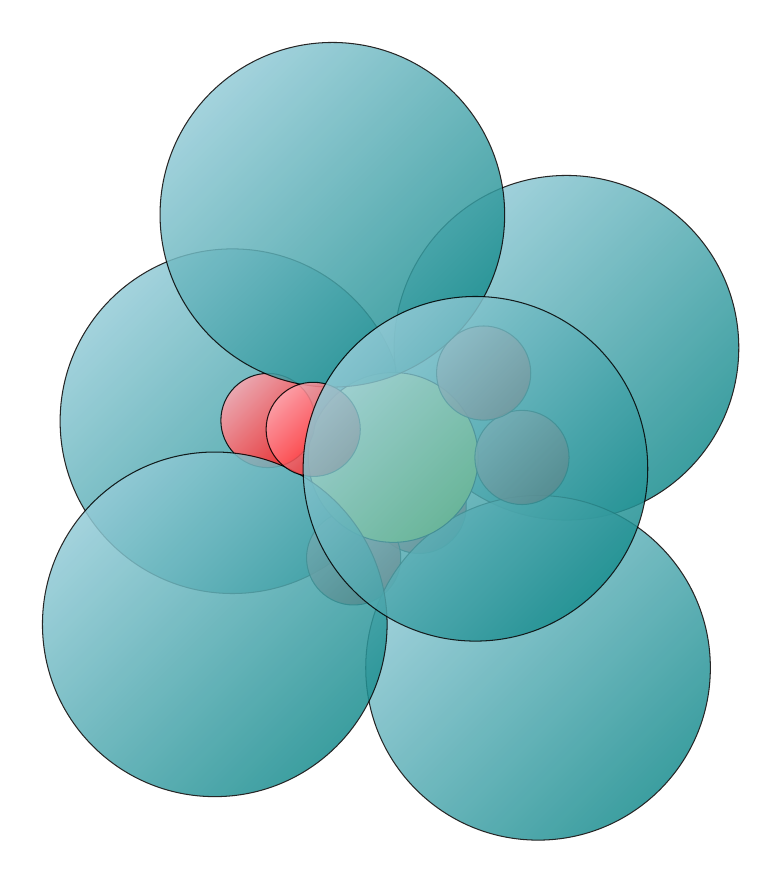}
\end{minipage}

\hspace{5mm}

\noindent
\begin{minipage}[t]{0.60\textwidth} \vspace{0pt}    
    \noindent\textbf{
        \noindent (20) 
    }\noindent\necklace{0}{1}{2221};
\necklace{0}{2}{2212211};
\necklace{0}{2}{2221211}

    \[
        (r_1, r_2) \approx (0.288, 1.033)
    \]

    \begin{dmath*}
    r_{1}^{8} - 4 r_{1}^{7} - 70 r_{1}^{6} - 212 r_{1}^{5} - 217 r_{1}^{4} - 76 r_{1}^{3} + 6 r_{1}^{2} + 8 r_{1} + 1
    \end{dmath*}

    \begin{dmath*}
    37 r_{2}^{8} - 102 r_{2}^{7} + 61 r_{2}^{6} + 48 r_{2}^{5} - 52 r_{2}^{4} - 16 r_{2}^{3} + 16 r_{2}^{2} + 8 r_{2} + 1    
    \end{dmath*}

\end{minipage} \begin{minipage}[t]{0.35\textwidth} \vspace{0pt}
   \includegraphics[width=4cm]{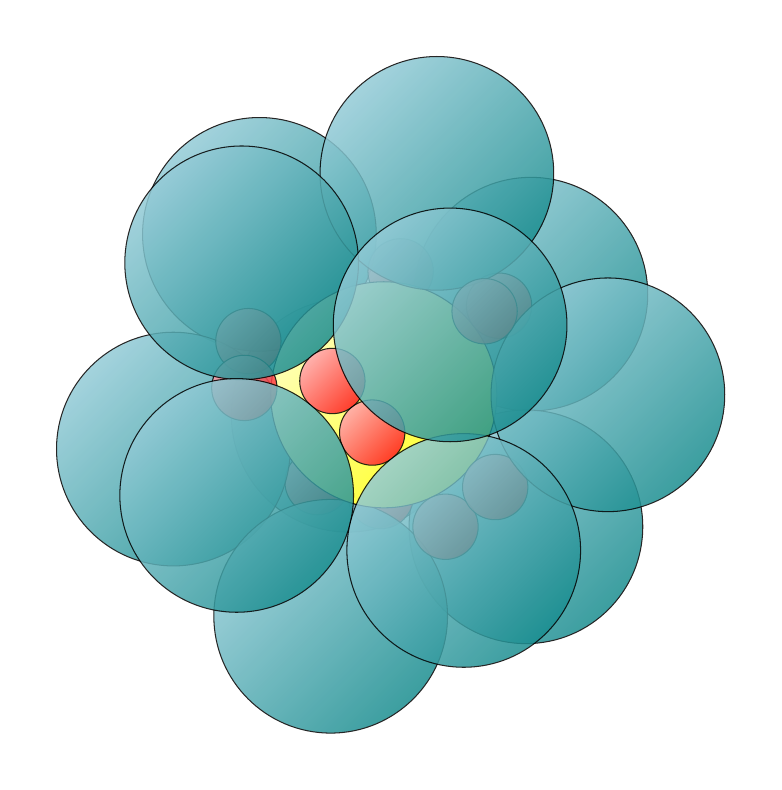}
\end{minipage}

\hspace{5mm}

\noindent
\begin{minipage}[t]{0.60\textwidth} \vspace{0pt}    
    \noindent\textbf{
        \noindent (21) 
    }\noindent\necklace{0}{1}{22211};
\necklace{0}{2}{221211}

    \[
        (r_1, r_2) \approx (0.465, 0.781)
    \]

    \begin{dmath*}
    3581 r_{1}^{22} - 9932 r_{1}^{21} - 61446 r_{1}^{20} + 111152 r_{1}^{19} + 552208 r_{1}^{18} - 239716 r_{1}^{17} - 2684838 r_{1}^{16} - 2138112 r_{1}^{15} + 4433136 r_{1}^{14} + 9568904 r_{1}^{13} + 5010912 r_{1}^{12} - 4127268 r_{1}^{11} - 7199879 r_{1}^{10} - 3409212 r_{1}^{9} + 555662 r_{1}^{8} + 1372800 r_{1}^{7} + 600696 r_{1}^{6} + 39792 r_{1}^{5} - 65480 r_{1}^{4} - 30400 r_{1}^{3} - 6448 r_{1}^{2} - 704 r_{1} - 32
    \end{dmath*}

    \begin{dmath*}
    256 r_{2}^{22} - 10272 r_{2}^{21} + 161236 r_{2}^{20} - 1222028 r_{2}^{19} + 4392369 r_{2}^{18} - 4995028 r_{2}^{17} - 7647390 r_{2}^{16} + 13650544 r_{2}^{15} + 7223766 r_{2}^{14} - 14580080 r_{2}^{13} - 6025260 r_{2}^{12} + 8739996 r_{2}^{11} + 4418825 r_{2}^{10} - 2805820 r_{2}^{9} - 2171890 r_{2}^{8} + 208768 r_{2}^{7} + 549976 r_{2}^{6} + 131232 r_{2}^{5} - 36080 r_{2}^{4} - 26560 r_{2}^{3} - 6256 r_{2}^{2} - 704 r_{2} - 32    
    \end{dmath*}

\end{minipage} \begin{minipage}[t]{0.35\textwidth} \vspace{0pt}
   \includegraphics[width=4cm]{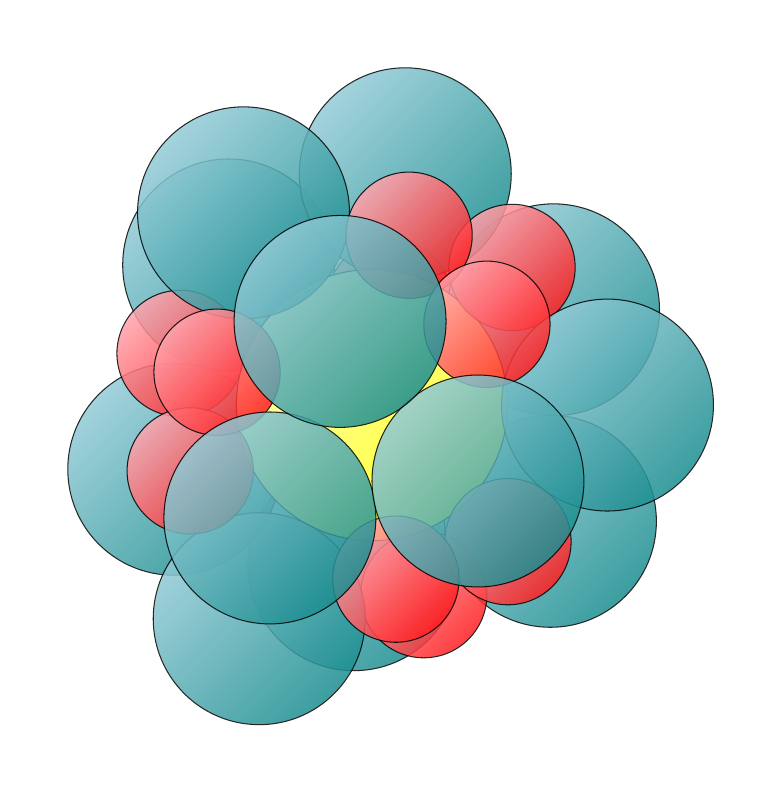}
\end{minipage}

\hspace{5mm}

\noindent
\begin{minipage}[t]{0.60\textwidth} \vspace{0pt}    
    \noindent\textbf{
        \noindent (22) 
    }\noindent\necklace{0}{1}{2222};
\necklace{0}{2}{212121}

    \[
        (r_1, r_2) \approx (0.577, 1.366)
    \]

    \begin{dmath*}
    3 r_{1}^{2} - 1
    \end{dmath*}

    \begin{dmath*}
    2 r_{2}^{2} - 2 r_{2} - 1    
    \end{dmath*}

\end{minipage} \begin{minipage}[t]{0.35\textwidth} \vspace{0pt}
   \includegraphics[width=4cm]{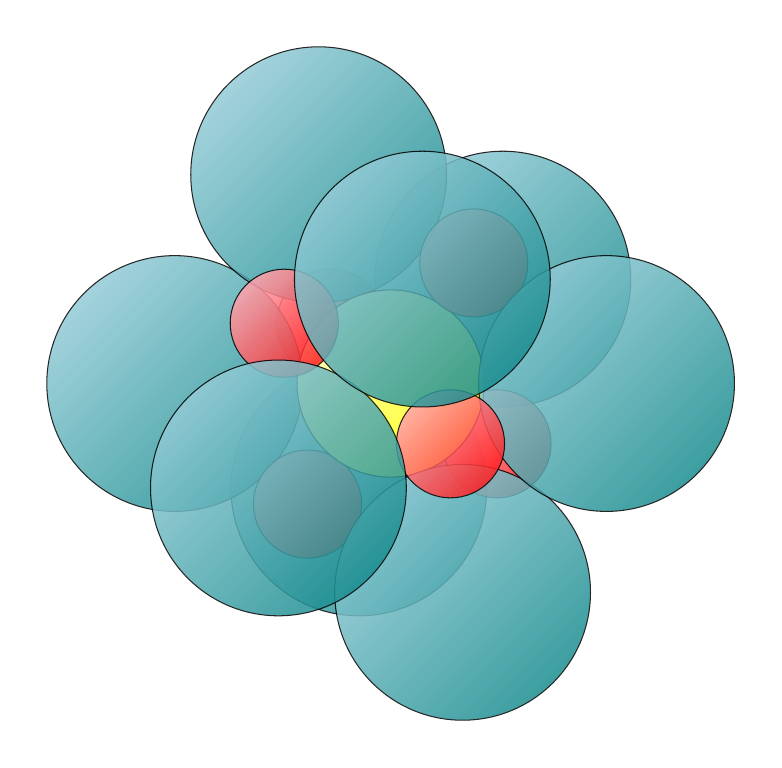}
\end{minipage}

\hspace{5mm}

\noindent
\begin{minipage}[t]{0.60\textwidth} \vspace{0pt}    
    \noindent\textbf{
        \noindent (23) 
    }\noindent\necklace{0}{1}{2222};
\necklace{0}{2}{22121}

    \[
        (r_1, r_2) \approx (0.905, 1.895)
    \]

    \begin{dmath*}
    r_{1}^{4} + 20 r_{1}^{3} - 2 r_{1}^{2} - 12 r_{1} - 3
    \end{dmath*}

    \begin{dmath*}
    4 r_{2}^{2} - 6 r_{2} - 3    
    \end{dmath*}

\end{minipage} \begin{minipage}[t]{0.35\textwidth} \vspace{0pt}
   \includegraphics[width=4cm]{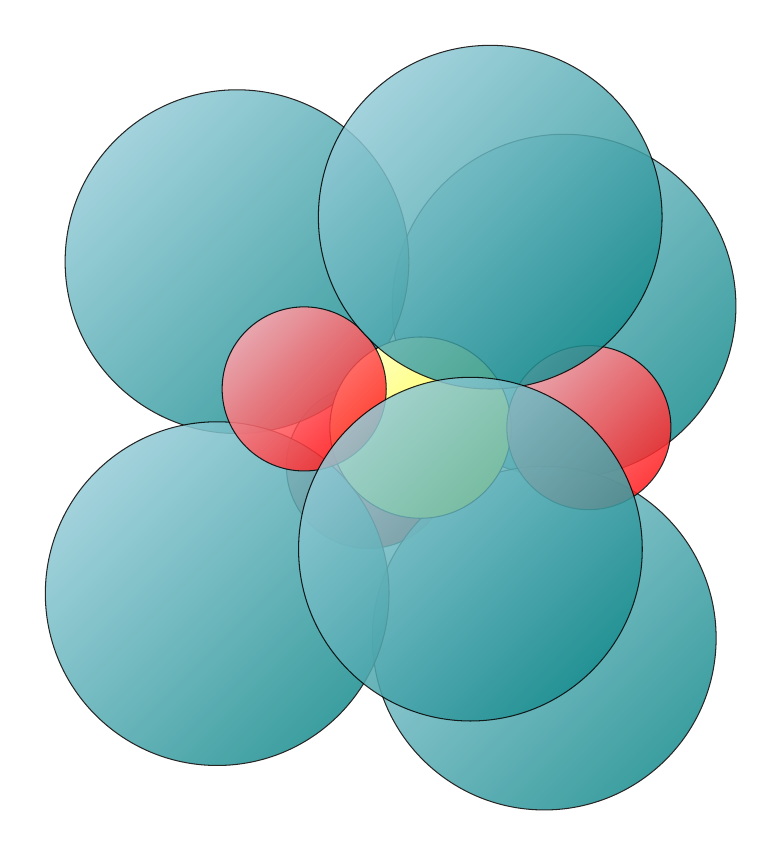}
\end{minipage}

\hspace{5mm}

\noindent
\begin{minipage}[t]{0.60\textwidth} \vspace{0pt}    
    \noindent\textbf{
        \noindent (24) 
    }\noindent\necklace{0}{1}{2222};
\necklace{0}{2}{2212121}

    \[
        (r_1, r_2) \approx (0.237, 0.556)
    \]

    \begin{dmath*}
    r_{1}^{4} - 16 r_{1}^{3} - 14 r_{1}^{2} + 1
    \end{dmath*}

    \begin{dmath*}
    8 r_{2}^{4} - 16 r_{2}^{3} - 4 r_{2}^{2} + 4 r_{2} + 1    
    \end{dmath*}

\end{minipage} \begin{minipage}[t]{0.35\textwidth} \vspace{0pt}
   \includegraphics[width=4cm]{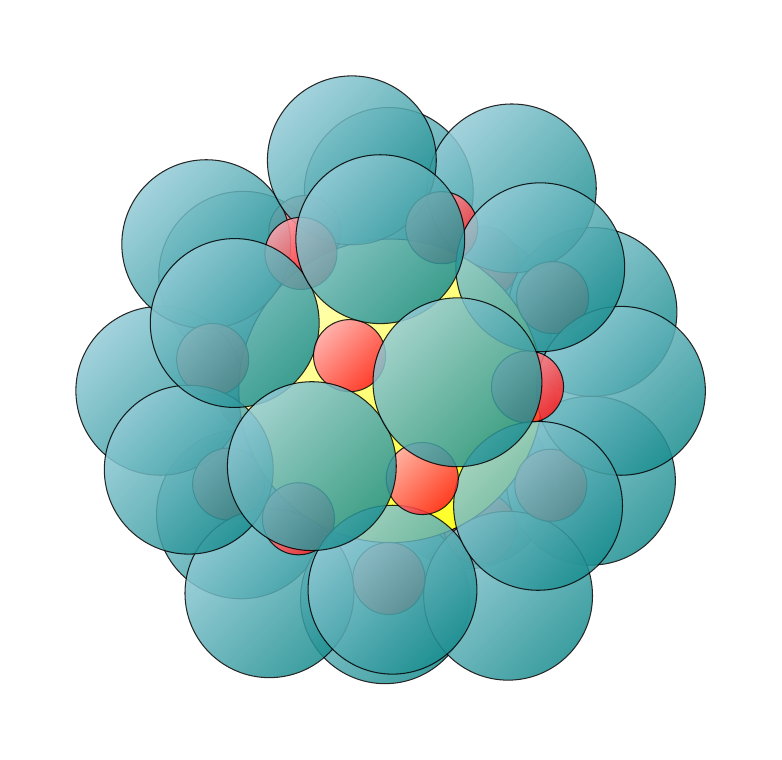}
\end{minipage}

\hspace{5mm}

\noindent
\begin{minipage}[t]{0.60\textwidth} \vspace{0pt}    
    \noindent\textbf{
        \noindent (25) 
    }\noindent\necklace{0}{1}{2222};
\necklace{0}{2}{221221};
\necklace{0}{2}{222121}

    \[
        (r_1, r_2) \approx (0.414, 1.000)
    \]

    \begin{dmath*}
    r_{1}^{2} + 2 r_{1} - 1
    \end{dmath*}

    \begin{dmath*}
    r_{2} - 1    
    \end{dmath*}

\end{minipage} \begin{minipage}[t]{0.35\textwidth} \vspace{0pt}
   \includegraphics[width=4cm]{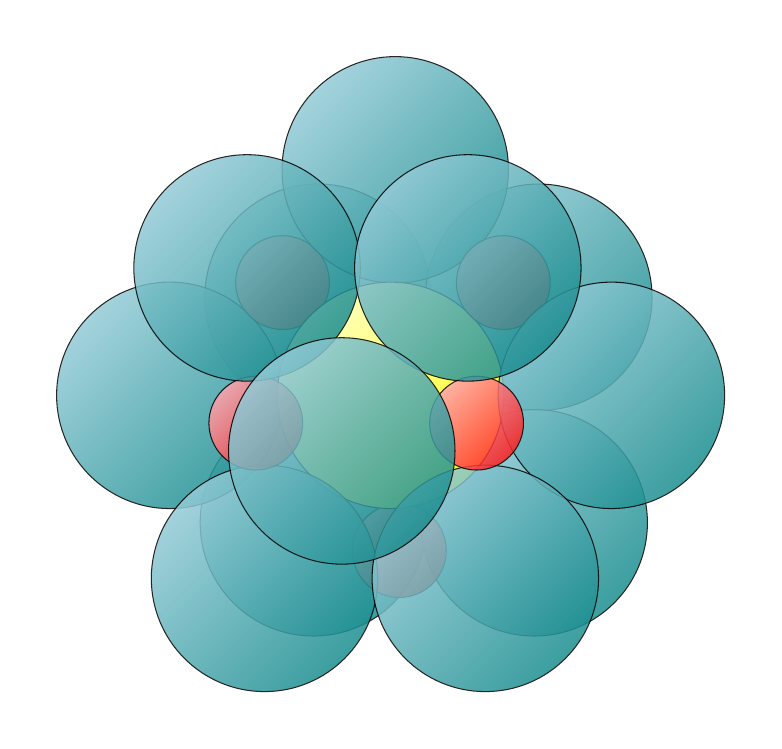}
\end{minipage}

\hspace{5mm}

\noindent
\begin{minipage}[t]{0.60\textwidth} \vspace{0pt}    
    \noindent\textbf{
        \noindent (26) 
    }\noindent\necklace{0}{1}{2222};
\necklace{0}{2}{22221}

    \[
        (r_1, r_2) \approx (0.672, 1.549)
    \]

    \begin{dmath*}
    7 r_{1}^{8} + 40 r_{1}^{7} + 28 r_{1}^{6} + 40 r_{1}^{5} + 34 r_{1}^{4} - 8 r_{1}^{3} - 20 r_{1}^{2} - 8 r_{1} - 1
    \end{dmath*}

    \begin{dmath*}
    r_{2}^{4} - 8 r_{2}^{3} + 4 r_{2}^{2} + 8 r_{2} + 2    
    \end{dmath*}

\end{minipage} \begin{minipage}[t]{0.35\textwidth} \vspace{0pt}
   \includegraphics[width=4cm]{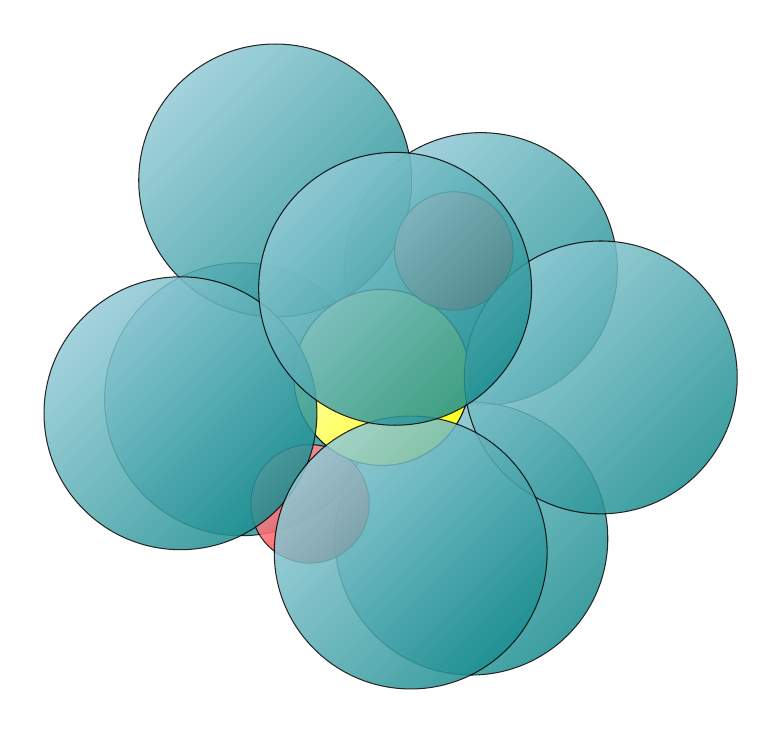}
\end{minipage}

\hspace{5mm}

\noindent
\begin{minipage}[t]{0.60\textwidth} \vspace{0pt}    
    \noindent\textbf{
        \noindent (27) 
    }\noindent\necklace{0}{1}{2222};
\necklace{0}{2}{222221}

    \[
        (r_1, r_2) \approx (0.251, 0.592)
    \]

    \begin{dmath*}
    239 r_{1}^{12} + 684 r_{1}^{11} + 1178 r_{1}^{10} + 2276 r_{1}^{9} + 2897 r_{1}^{8} + 2840 r_{1}^{7} + 2668 r_{1}^{6} + 1832 r_{1}^{5} + 673 r_{1}^{4} + 60 r_{1}^{3} - 38 r_{1}^{2} - 12 r_{1} - 1
    \end{dmath*}

    \begin{dmath*}
    r_{2}^{6} - 10 r_{2}^{5} + 23 r_{2}^{4} + 20 r_{2}^{3} - 5 r_{2}^{2} - 6 r_{2} - 1    
    \end{dmath*}

\end{minipage} \begin{minipage}[t]{0.35\textwidth} \vspace{0pt}
   \includegraphics[width=4cm]{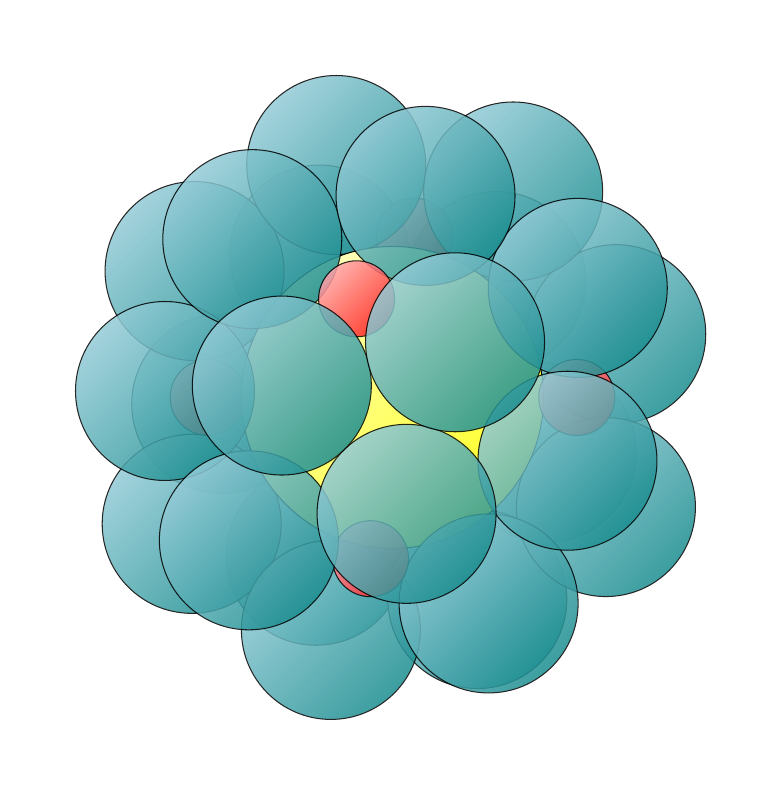}
\end{minipage}

\hspace{5mm}

\noindent
\begin{minipage}[t]{0.60\textwidth} \vspace{0pt}    
    \noindent\textbf{
        \noindent (28) 
    }\noindent\necklace{0}{1}{22222};
\necklace{0}{2}{212121}

    \[
        (r_1, r_2) \approx (0.403, 0.554)
    \]

    \begin{dmath*}
    59 r_{1}^{8} + 22 r_{1}^{7} - 213 r_{1}^{6} - 96 r_{1}^{5} + 120 r_{1}^{4} + 64 r_{1}^{3} - 8 r_{1}^{2} - 8 r_{1} - 1
    \end{dmath*}

    \begin{dmath*}
    r_{2}^{4} - 14 r_{2}^{3} - 3 r_{2}^{2} + 4 r_{2} + 1    
    \end{dmath*}

\end{minipage} \begin{minipage}[t]{0.35\textwidth} \vspace{0pt}
   \includegraphics[width=4cm]{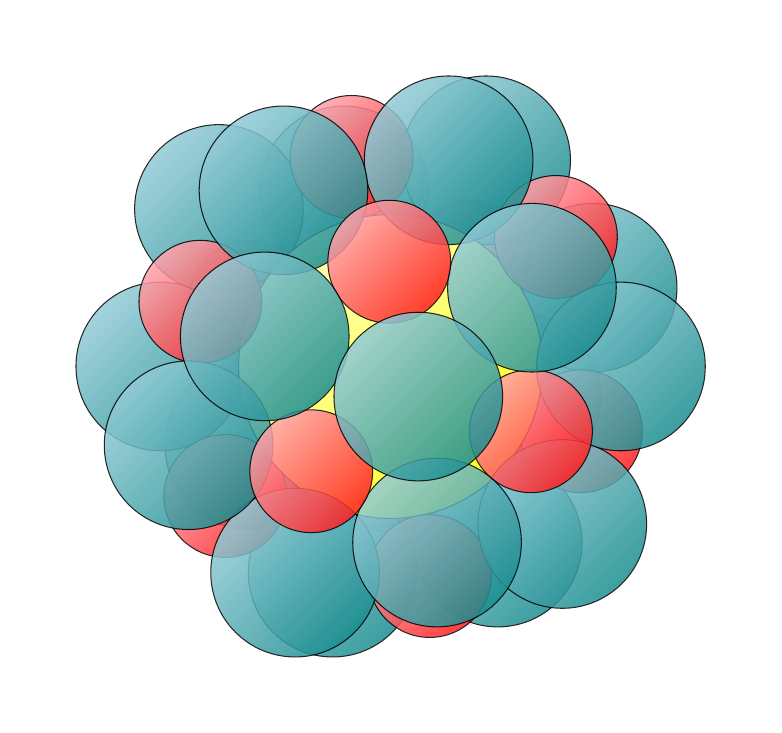}
\end{minipage}

\hspace{5mm}

\noindent
\begin{minipage}[t]{0.60\textwidth} \vspace{0pt}    
    \noindent\textbf{
        \noindent (29) 
    }\noindent\necklace{0}{1}{22222};
\necklace{0}{2}{221221}

    \[
        (r_1, r_2) \approx (0.318, 0.447)
    \]

    \begin{dmath*}
    19 r_{1}^{4} + 36 r_{1}^{3} + 9 r_{1}^{2} - 4 r_{1} - 1
    \end{dmath*}

    \begin{dmath*}
    5 r_{2}^{2} - 1    
    \end{dmath*}

\end{minipage} \begin{minipage}[t]{0.35\textwidth} \vspace{0pt}
   \includegraphics[width=4cm]{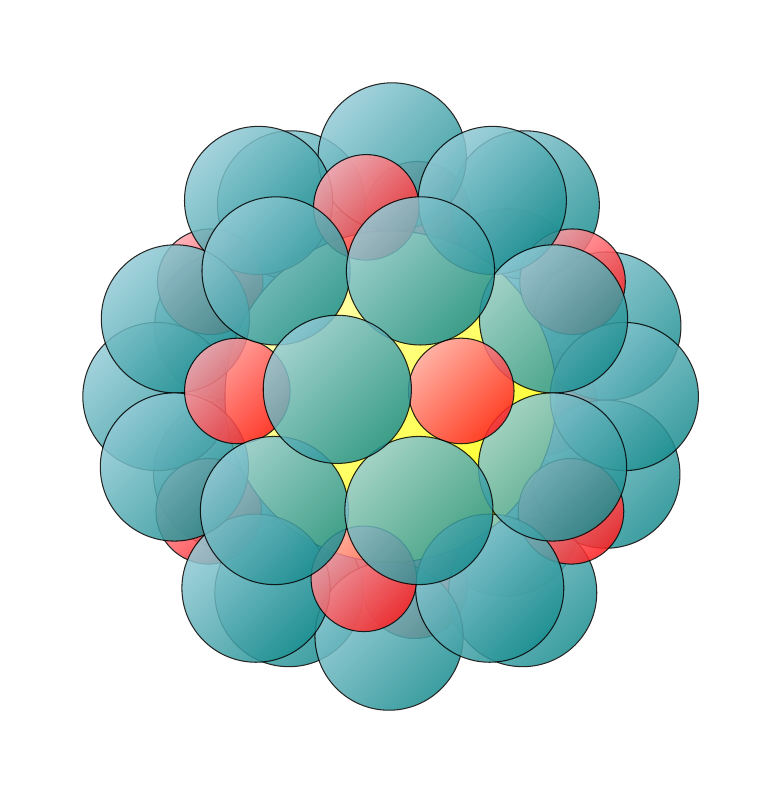}
\end{minipage}

\hspace{5mm}

\noindent
\begin{minipage}[t]{0.60\textwidth} \vspace{0pt}    
    \noindent\textbf{
        \noindent (30) 
    }\noindent\necklace{0}{1}{22222};
\necklace{0}{2}{222221}

    \[
        (r_1, r_2) \approx (0.212, 0.301)
    \]

    \begin{dmath*}
    6479 r_{1}^{24} + 381546 r_{1}^{23} + 1543709 r_{1}^{22} - 2733224 r_{1}^{21} - 30595326 r_{1}^{20} - 100099884 r_{1}^{19} - 204383266 r_{1}^{18} - 293842064 r_{1}^{17} - 303990446 r_{1}^{16} - 217207604 r_{1}^{15} - 84894026 r_{1}^{14} + 15536176 r_{1}^{13} + 45814554 r_{1}^{12} + 28542506 r_{1}^{11} + 6312019 r_{1}^{10} - 2199984 r_{1}^{9} - 1889791 r_{1}^{8} - 437364 r_{1}^{7} + 30954 r_{1}^{6} + 36796 r_{1}^{5} + 6429 r_{1}^{4} - 154 r_{1}^{3} - 191 r_{1}^{2} - 24 r_{1} - 1
    \end{dmath*}

    \begin{dmath*}
    r_{2}^{12} - 20 r_{2}^{11} + 146 r_{2}^{10} - 428 r_{2}^{9} + 187 r_{2}^{8} + 938 r_{2}^{7} - 153 r_{2}^{6} - 976 r_{2}^{5} - 536 r_{2}^{4} - 50 r_{2}^{3} + 39 r_{2}^{2} + 12 r_{2} + 1    
    \end{dmath*}

\end{minipage} \begin{minipage}[t]{0.35\textwidth} \vspace{0pt}
   \includegraphics[width=4cm]{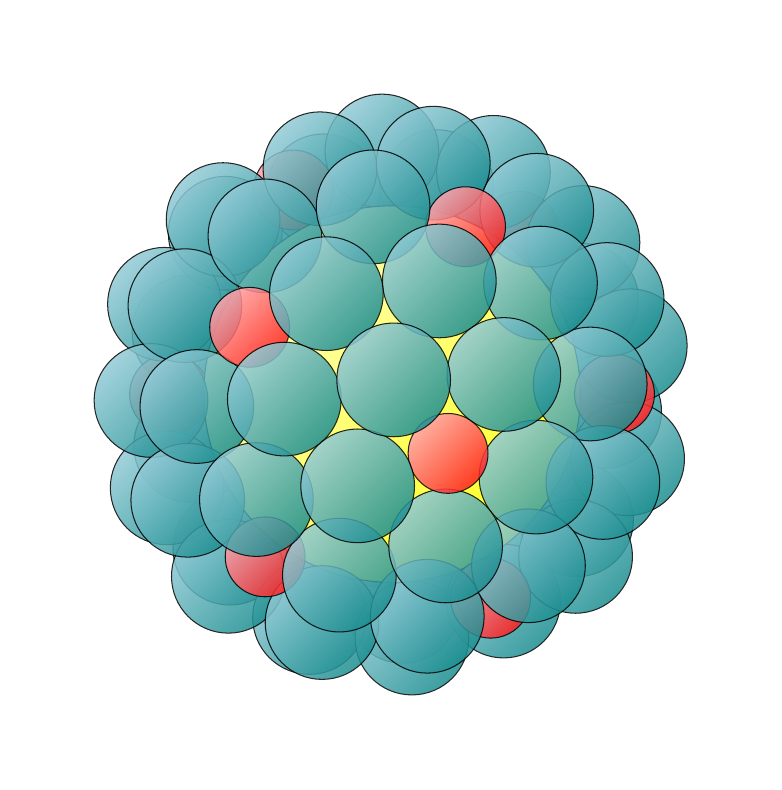}
\end{minipage}

\hspace{5mm}

\section{Data availability statement}

Two datasets were generated in support of the results in this paper.
These datasets are described in Section~\ref{sec:Computing-all-raspberries}
and are publicly available \cite{raspberrydata}.

\bibliographystyle{amsalpha}
\bibliography{bibliography}

@article{Kennedy2006,
  author   = {Kennedy, T.},
  title    = {Compact packings of the plane with two sizes of discs},
  journal  = {Discrete Comput. Geom.},
  fjournal = {Discrete \& Computational Geometry. An International Journal
              of Mathematics and Computer Science},
  volume   = {35},
  year     = {2006},
  number   = {2},
  pages    = {255--267},
  mrnumber = {2195054},
  doi      = {10.1007/s00454-005-1172-4}
}

@article{Messerschmidt2020,
  author   = {Messerschmidt, M.},
  title    = {On compact packings of the plane with circles of three radii},
  journal  = {Comput. Geom.},
  fjournal = {Computational Geometry. Theory and Applications},
  volume   = {86},
  year     = {2020},
  doi      = {10.1016/j.comgeo.2019.05.002}
}

@article{Fernique2021,
  author   = {Fernique, T. and Hashemi, A. and Sizova, O.},
  title    = {Compact {P}ackings of the {P}lane with {T}hree {S}izes of
              {D}iscs},
  journal  = {Discrete Comput. Geom.},
  fjournal = {Discrete \& Computational Geometry. An International Journal
              of Mathematics and Computer Science},
  volume   = {66},
  year     = {2021},
  number   = {2},
  pages    = {613--635},
  issn     = {0179-5376},
  mrnumber = {4292755},
  doi      = {10.1007/s00454-019-00166-y}
}

@article{FerniqueTwoSpheres2021,
  author   = {Fernique, T.},
  title    = {Compact packings of space with two sizes of spheres},
  journal  = {Discrete Comput. Geom.},
  fjournal = {Discrete \& Computational Geometry. An International Journal
              of Mathematics and Computer Science},
  volume   = {65},
  year     = {2021},
  number   = {4},
  pages    = {1287--1295},
  mrclass  = {52C17},
  mrnumber = {4249904},
  doi      = {10.1007/s00454-019-00140-8}
}

@article{FerniqueThreeSizes,
  author        = {Fernique, T.},
  title         = {Compact packings of space with three sizes of spheres},
  journal       = {https://arxiv.org/abs/1912.02293},
  year          = {2019},
  archiveprefix = {arXiv},
  eprint        = {1912.02293}
}

@article{Messerschmidt2d,
  author   = {Messerschmidt, M.},
  title    = {The number of configurations of radii that can occur in compact packings of the plane with discs of n sizes is finite},
  journal  = {Discrete Comput. Geom.},
  fjournal = {Discrete \& Computational Geometry. An International Journal
              of Mathematics and Computer Science},
  year     = {2023},
  doi      = {10.1007/s00454-022-00471-z}
}

@article{WinterPolytopes,
  author  = {Winter, M.},
  title   = {Rigidity, Tensegrity, and Reconstruction of Polytopes Under Metric Constraints},
  journal = {International Mathematics Research Notices},
  volume  = {2024},
  number  = {9},
  pages   = {7721-7747},
  year    = {2023},
  doi     = {10.1093/imrn/rnad298}
}

@article{ChemistryFernique,
  author  = {Chinaud-Chaix, C. and Marchenko, N. and Fernique, T. and Tricard, S.},
  title   = {Do chemists control plane packing{,} i.e. two-dimensional self-assembly{,} at all scales?},
  journal = {New J. Chem.},
  year    = {2023},
  volume  = {47},
  issue   = {15},
  pages   = {7014-7025},
  doi     = {10.1039/D3NJ00208J},
  url     = {http://dx.doi.org/10.1039/D3NJ00208J}
}

@article{ChemistryPaik,
  author  = {Paik, T. and Diroll, B.T. and Kagan, C.R. and Murray, C.B.},
  title   = {Binary and Ternary Superlattices Self-Assembled from Colloidal Nanodisks and Nanorods},
  journal = {Journal of the American Chemical Society},
  volume  = {137},
  number  = {20},
  pages   = {6662-6669},
  year    = {2015},
  doi     = {10.1021/jacs.5b03234}
}

@article{MesserschmidtKikianty2024,
  author   = {M. Messerschmidt and E. Kikianty},
  title    = {On Compact Packings of Euclidean Space with Spheres of Finitely Many Sizes},
  journal  = {Discrete Comput. Geom.},
  fjournal = {Discrete \& Computational Geometry. An International Journal
              of Mathematics and Computer Science},
  year     = {2024},
  doi      = {10.1007/s00454-024-00628-y}
}

@article{MathewsZymaris,
  author   = {Mathews, D.V. and Zymaris, O.},
  title    = {Spinors and the {D}escartes circle theorem},
  journal  = {J. Geom. Phys.},
  fjournal = {Journal of Geometry and Physics},
  volume   = {212},
  year     = {2025},
  pages    = {Paper No. 105458, 14},
  doi      = {10.1016/j.geomphys.2025.105458}
}

@book{AdamsLoustaunau,
  author    = {Adams, W.W. and Loustaunau, P.},
  title     = {An introduction to {G}r\"{o}bner bases},
  publisher = {American Mathematical Society, Providence, RI},
  year      = {1994},
  doi       = {10.1090/gsm/003}
}

@article{sympy,
  title   = {SymPy: symbolic computing in Python},
  author  = {Meurer, A. and Smith, C. P. and Paprocki, M. and \v{C}ert\'{i}k, O. and Kirpichev, S. B. and Rocklin, M. and Kumar, A. and Ivanov, S. and Moore, J. K. and Singh, S. and Rathnayake, T. and Vig, S. and Granger, B. E. and Muller, R. P. and Bonazzi, F. and Gupta, H. and Vats, S. and Johansson, F. and Pedregosa, F. and Curry, M. J. and Terrel, A. R. and Rou\v{c}ka, \v{S}. and Saboo, A. and Fernando, I. and Kulal, S. and Cimrman, R. and Scopatz, A.},
  year    = 2017,
  volume  = 3,
  pages   = {e103},
  journal = {PeerJ Computer Science},
  doi     = {10.7717/peerj-cs.103}
}

@misc{SINGULAR,
  title        = {{\sc Singular} {4-4-0} --- {A} computer algebra system for polynomial computations},
  author       = {Decker, W. and Greuel, G-M. and Pfister, G. and Sch\"onemann, H.},
  year         = {2024},
  howpublished = {\url{http://www.singular.uni-kl.de}}
}

@misc{mpmath,
  author       = {F. Johansson},
  title        = {mpmath: A Python library for arbitrary-precision floating-point arithmetic},
  year         = {2007},
  url          = {https://mpmath.org/},
  howpublished = {\url{https://mpmath.org/}}
}

@misc{raspberrydata,
  author       = {M. Messerschmidt},
  title        = {Data for {R}aspberries with at most two sizes of berry},
  year         = {2025},
  howpublished = {\url{https://github.com/miekmesserschmidt/data-raspberries-with-two-sizes-of-berry}}
}
 \end{document}